%% file: HGHT.tex
\documentclass[twoside]{irmaems}
\usepackage{amssymb} 
\usepackage{amsmath} 
\usepackage{latexsym}

\usepackage{tikz}
\usetikzlibrary{calc}
\usepackage{pgf,tikz}
\usetikzlibrary{arrows}

\usepackage{tikz}
\usepackage{ color, amsmath, amssymb, amsfonts, amscd, graphicx,latexsym}
\usepackage{caption, subcaption, amsfonts}
\usepackage[pdf]{pstricks}
\usepackage{pst-grad} 
\newcommand{\ITEM}{$\blacktriangleright$}

\newcommand{\gal}{G_{\mathbb Q}}

\newcommand{\Q}{\mathbb Q}
\newcommand{\Z}{\mathbb Z}

\newcommand{\CC}{\mathbf C}
\newcommand{\p}{\mathbb P}

\newcommand{\psl}{\mathrm{PSL}_2(\Z)}
\def\modorb{\mbox{$\circ\hspace{-1.5mm}-\!\!\!-\hspace{-1.5mm}\bullet$}}

\def\ztwo{\Z\!/\!2\Z}

\def\zfour{\Z\!/\!4\Z}

\renewcommand\theequation{\arabic{section}.\arabic{equation}}

\theoremstyle{definition} 

 \newtheorem{definition}{Definition}[section]

\theoremstyle{plain}      

 \newtheorem{theorem}[definition]{Theorem}

\newtheorem*{theorem*}{Theorem}

\newtheorem*{definition*}{Definition}

\begin{document}

\title{Hypergeometric Galois Actions}

\author{Muhammed Uluda\u{g}\hspace{.1mm} \thanks{
Work supported by T\"{U}B\.{I}TAK Grant No.~110T690 and a GSU research grant.} \, and  \.{I}smail Sa\u{g}lam \thanks{Work partially
supported by T\"{U}B\.{I}TAK Grant No.~110T690.}}

\address{
Galatasaray University\\
\c{C}{\i}ra\u{g}an Cad. No: 36, 34349, Ortak\"{o}y, \.{I}stanbul, Turkey \\
email:\,\tt{muhammed.uludag@gmail.com}
\\[4pt]
Galatasaray University\\
\c{C}{\i}ra\u{g}an Cad. No: 36, 34349, Ortak\"{o}y, \.{I}stanbul, Turkey \\
email:\,\tt{isaglamtrfr@gmail.com}}

\maketitle

\begin{abstract} 
We outline a project to study the Galois action on a class of modular graphs (special type of dessins) 
which arise as the dual graphs of the sphere triangulations of non-negative curvature, classified by Thurston.
Because of their connections to hypergeometric functions, there is a hope that these graphs will render
themselves to explicit calculation for a study of Galois action on them, unlike the case of a general dessin.
\end{abstract}


\begin{keywords}
Sphere triangulation, hypergeometric functions, dessins, Belyi maps, modular graphs, trivalent ribbon graphs, Galois actions, cone metric, flat structure, euclidean structure, ball quotient, branched covering of the sphere, complex hyperbolic space.
\end{keywords}

\tableofcontents   

\section{Introduction}
How to get useful information about the absolute Galois group from dessins?  In order to reply to this question, i.e. to compute the Galois action on a dessin, we need to compute its Belyi map. This problem is algorithmically solvable, but often returns some complicated expressions which are hard to treat in a systematic manner in the full generality of the problem. On the other hand, even if we are able to compute the Galois action on an individual dessin, this is just a finite action of $\gal$ and cannot yield information about its profinite structure.

We are thus led to seek some special infinite families of dessins which can be studied in a systematic manner. We may reformulate this problem in terms of the coverings 
$$
X\rightarrow \mathbb P^1(\mathbf C)\setminus\{0,1,\infty\}
$$
of the thrice-punctured sphere. As is well-known, these coverings correspond in a 1-1 manner to dessins. In terms of coverings, we are interested in infinite ``systems" of essentially non-abelian coverings.

The thrice-punctured sphere has the standard ideal triangulation which consists of two triangles with vertices at $0,1,\infty$. 
Lifting this triangulation via the covering map, we obtain a triangulation of the covering surface $X$. The idea of the present paper is to impose a cone metric on $X$ by declaring these triangles to be congruent euclidean (flat) equilateral  triangles. We are interested in the case where $X$ is a punctured sphere (the corresponding dessin being a dessin on a punctured sphere).

Thus we have a punctured sphere with an ideal triangulation, and we are led to the question: is it possible to understand sphere triangulations in a systematic manner? It turns out that, if we impose a certain ``non-negative curvature" condition on the induced cone metric, then answer to this question is very positive. These triangulations are parametrized by the points lying inside a cone\footnote{Beware  the use of the word ``{cone}" in two distinct senses.} in a certain 20-dimensional integral lattice modulo some automorphism group of the lattice. They can be explicitly constructed by  cut-and-glue operations. 

Not every triangulation comes from a covering, but there is a remedy for this problem, by considering the graph dual to the triangulation. We start the Section 2 at this point, and show that a triangulation is nothing but a covering of the modular curve. Section 3 introduces the metric point of view and provides the first contact with Thurston's classification. In addition, we point out to some amusing connections with chemistry and the genus-0 phenomenon of moonshine. In Section 4 we come back to the covering interpretation of triangulations and present a simple application of the Riemann-Hurwitz formula.
As a result we rediscover the famous list of integer tuples (Appendix 1) due to Terada, Deligne\&Mostow, reproduced in an alternative way by Thurston. We speculate on the existence of other types of classifiable branching problems and perform some numerology. Results are given in Appendix 2-3. Section 5 is devoted to an exposition of Thurston's theory and also provides a contact with hypergeometric functions. The section ends with a series of problems related to arithmetic aspects. Section 6 is an exposition of a chapter of \.Ismail Sa\u glam's thesis \cite{ismailthesis} and gives a case study of the simplest ``system" of triangulations. In Section 7 we shortly explain how one can go beyond Thurston's classification.

As its name suggests, this quest aspires to be a continuation of the ``Geometric Galois Actions" initiative of Schneps and Lochak \cite{gga1}, \cite{gga2}, \cite{gga3} from the 90's. The paper by Zvonkine and Magot \cite{bafas} is another precursor of our approach in that it studies the Belyi maps related to some Archimedean polyhedra, a few being related to the triangulations of non-negative curvature. 
To our knowledge, besides our work
\cite{ayberkdessin3}, \cite{ayberkmakale}, \cite{ismail}, \cite{ismailthesis} there are no other 
attempts to realize Grothendieck's dream in the hypergeometric context.

\section{Category of coverings of the modular curve}
For more details on this section, see \cite{A Panorama of the Fundamental Group of the Modular Curve Uludag-Zeytin}.
Our aim is here to establish an equivalence between triangulations of surfaces and the bipartite dual graphs, constructed by putting a vertex of type $\bullet$ at the center of each triangle, connecting these vertices via edges and putting a vertex of type $\circ$ whenever this edge meets an arc of the triangulation\footnote{We require that an edge and an arc meets always transversally and at most at one point. Also note that we are interested in combinatorial types (i.e. homeomorphism classes) of triangulations and graphs.}.
We call these graphs {\it modular graphs}, including the duals of degenerate triangulations. If the triangulation is finite and consists of $n$ non-degenerate triangle, then its dual modular has $3n$ edges.

Modular graphs constitute a special class of dessins. Just as dessins classify the conjugacy classes of subgroups of the thrice-punctured sphere, modular graphs classify the conjugacy classes subgroups of the modular group. This correspondence extends to a correspondence between modular graphs with a chosen edge and subgroups (i.e. not only conjugacy classes) of the modular group.  Denote by ${\rm {\bf FSub} }(\psl)$ the category of all finite-index subgroups of $\psl$, with inclusions as morphisms.
Our claim is that (pointed) modular graphs constitute a category with coverings as morphisms, and the pointed former category is equivalent to the category ${\rm {\bf FSub} }(\psl)$.

Consider the arc connecting the two elliptic points on the boundary of the standard fundamental domain of the $\psl$ action on $\mathbb H$. Then the $\psl$-orbit of this arc is a tree $\mathcal F$, called the {\it Farey tree}\index{Farey tree}. This tree admits a $\psl$-action by definition, and the quotient graphs by subgroups of finite or infinite index $G<\psl$ gives precisely the {\it modular graphs}\index{modular graph} \cite{ayberkdessin1} introduced above as duals of triangulations. In particular, the quotient \emph{orbi-graph} $\mathcal F/\psl$ is an arc connecting the two orbifold points of the modular orbifold $\mathbb H/\psl$. We call this the {\it modular arc}\index{modular arc} and denote it by $\modorb$. Its (pointed) covering category is defined 
respectively by ${\rm {\bf FCov} } ^*(\modorb)$ and  ${\rm {\bf FCov} }(\modorb)$, and consists precisely of modular graphs, i.e. duals graphs of triangulations including degenerate ones. The claimed equivalence follows. 

The quotient of the upper half plane under the $\psl$ action is called the {\it modular orbifold}\footnote{Also known by the names {\it modular curve} or {\it modular surface}.} and denoted $\mathcal M$. It can be identified with the sphere with a puncture at infinity and  with two orbifold points 0 and 1 with $\Z/2\Z$ and $\Z/3\Z$-inertia respectively. 
The fundamental group of the modular orbifold is $\psl$.
By the usual correspondence from topology, 
its pointed covering category ${\rm {\bf FCov}^* }(\mathcal M)$ is arrow-reversing equivalent to the category 
${\rm {\bf FSub} }(\psl)$. Since this latter is precisely the category of modular graphs, we see that the modular graphs classify the coverings of the modular orbifold.

Since the modular graphs are dual to surface triangulations, we see that a non-degenerate surface triangulation with $n$ triangles is nothing but a degree $3n$-covering of the modular orbifold.

As the simplest instance of this correspondence, recall that the congruence modular group $\Gamma(2)<\psl$ acts on $\mathbb H$ freely, the quotient being the thrice-punctured sphere $\mathbb P^1(\mathbf C)\setminus\{0,1,\infty\}$. This sphere admits a unique ideal triangulation with two triangles. The dual modular graph has six edges, two type-$\bullet$ and three type-$\circ$  vertices.
So we rediscover the well-known fact that, $\mathbb P^1(\mathbf C)\setminus\{0,1,\infty\}$ is a degree-6 covering of the modular orbifold. 

\section{Clash of Geometrizations}
Until now, we had algebra, arithmetic and combinatorics in the picture, but we have not made an essential use of a metric.

Being a quotient of the upper half plane under the action of a subgroup  $G<\mathrm{PSL}(2, \mathbb{Z})$ which preserves the hyperbolic metric, every 
surface $\mathbb H/G$ carries a canonical hyperbolic metric. This is a punctured surface and the metric becomes infinite at the cusps.
We have seen that the covering $\mathbb H/G\to \mathcal M$ is also determined by  the combinatorial class  of an ideal triangulation, (including degenerate ones) with vertices at the cusps.  
Now we introduce a flat metric on the modular orbifold $\mathcal M$, as follows. First put the flat metric on the canonical ideal triangulation of  $\mathbb P^1(\mathbf C)\setminus\{0,1,\infty\}$ by identifying its triangles by a equilateral euclidean triangle.
This metric also admits a $\Sigma_3$ symmetry and defines a metric on the quotient surface $\mathcal M$.
This metric lifts to every covering of $\mathcal M$ and this way every $\mathbb H/G$ becomes  $\diamond$ an equilateral-triangulated surface 
(for degenerate triangulations one must modify this claim a bit). For example, the thrice-punctured sphere becomes equilaterally triangulated with two equilateral triangles with vertices at the cusps 0,1,$\infty$.

There is an abrupt change of geometry in the above paragraph which precisely occurs at the
$\diamond$ sign: every surface $\mathbb H/G$ has been endowed with a Euclidean structure.
 Is this a natural structure? Yes, if you think that it is natural to identify the modular tile with the equilateral triangle modulo $\Sigma_3$. But somebody else may find it natural to identify it with a spherical triangle, see \cite{fengluo}. 
 
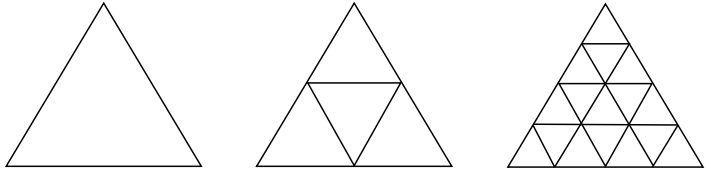
\begin{figure}[ht]\begin{center}
\scalebox{.5} 
{
\begin{pspicture}(0,-2.21)(18.588829,2.191018)
\pstriangle[linewidth=0.04,dimen=outer](2.63,-2.19)(5.26,4.4)
\pstriangle[linewidth=0.04,dimen=outer](9.29,-2.19)(5.26,4.4)
\rput{-180.0}(18.58,-2.12){\pstriangle[linewidth=0.04,dimen=outer](9.29,-2.19)(2.54,2.26)}
\pstriangle[linewidth=0.04,dimen=outer](15.97,-2.21)(5.26,4.4)
\rput{-180.0}(31.94,-2.16){\pstriangle[linewidth=0.04,dimen=outer](15.97,-2.21)(2.54,2.26)}
\psline[linewidth=0.04cm](14.62,-2.19)(16.62,1.11)
\psline[linewidth=0.04cm](15.34,1.09)(17.24,-2.15)
\psline[linewidth=0.04cm](14.06,-1.05)(17.9,-1.07)
\psline[linewidth=0.04cm](14.62,-2.19)(14.04,-1.07)
\psline[linewidth=0.04cm](17.22,-2.19)(17.9,-1.07)
\psline[linewidth=0.04cm](15.34,1.09)(16.62,1.09)
\end{pspicture} 
}
\caption{The simplest hypergeometric sphere triangulations -- the sphere is obtained by gluing two copies of the triangles along their boundaries.}\end{center}
\end{figure}
 
We {\it do} think that this structure is useful from the point of view of arithmetic. For example, there is a natural operation on the set of equilateral triangulations, i.e. the simultaneous subdivision of all its triangles, see the picture below. Note that this operation adds new vertices (cusps) to the triangulation. Although very neatly organized with respect to each other, these triangulations do not constitute a chain of coverings inside $\mathbf{FCov}(\mathcal M)$. Nevertheless, thanks to their connections with elliptic curves, we have succeeded in determining their Belyi maps in terms of the Weierstrass $\mathcal P$-function \cite{ayberkdessin1} (the same for the quadrangulations below, \cite{ayberkdessin2}). 
Hence, {\it this is a new kind of natural structure inside the category $\mathbf{FCov}(\mathcal M)$,} which has its origins in geometry; or rather hypergeometry, as we shall see.

\begin{figure}[h]\begin{center}
\scalebox{.5} 
{
\begin{pspicture}(0,-2.54)(18.72,2.56)
\psframe[linewidth=0.04,dimen=outer](5.02,2.52)(0.0,-2.5)
\psframe[linewidth=0.04,dimen=outer](11.84,2.54)(6.82,-2.48)
\psline[linewidth=0.04cm](9.34,2.52)(9.36,-2.52)
\psline[linewidth=0.04cm](6.84,0.02)(11.84,0.02)
\psframe[linewidth=0.04,dimen=outer](18.68,2.54)(13.66,-2.48)
\psline[linewidth=0.04cm](16.18,2.52)(16.2,-2.52)
\psline[linewidth=0.04cm](13.68,0.02)(18.68,0.02)
\psline[linewidth=0.04cm](14.94,2.54)(14.96,-2.48)
\psline[linewidth=0.04cm](17.46,2.52)(17.48,-2.48)
\psline[linewidth=0.04cm](13.68,1.28)(18.68,1.26)
\psline[linewidth=0.04cm](13.66,-1.24)(18.7,-1.26)
\end{pspicture} 
}
\caption{The simplest hypergeometric sphere quadrangulations -- the sphere is obtained by gluing two copies of the quadrangles along their boundaries.}\end{center}
\end{figure}
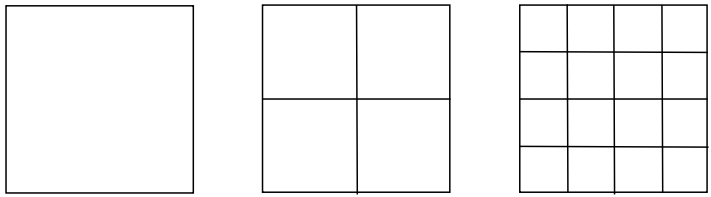

\begin{tabular}{|c|c|c|c|c|}
\hline
\begin{tabular}{c}cone\\picture\end{tabular}&
\begin{tabular}{c}$d$\\(vertex degree)\end{tabular}&
$6-d$&
\begin{tabular}{c}$\kappa$\\(curvature)\end{tabular}&
\begin{tabular}{c}$\theta=2\pi-\kappa$\\(cone angle)\end{tabular}\\
\hline
\includegraphics[scale=.1]{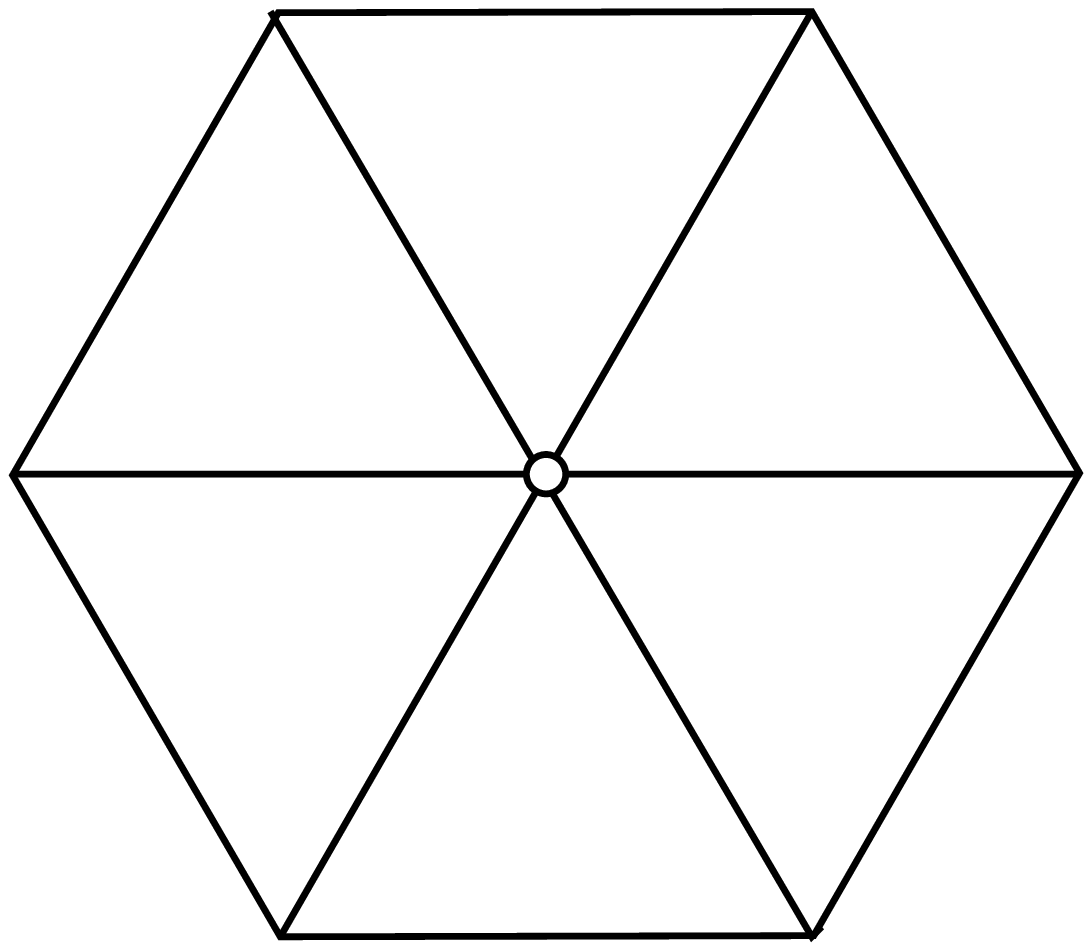}&
6&
0&
$\kappa=0$&
$2\pi$\\ 
\hline
\includegraphics[scale=.1]{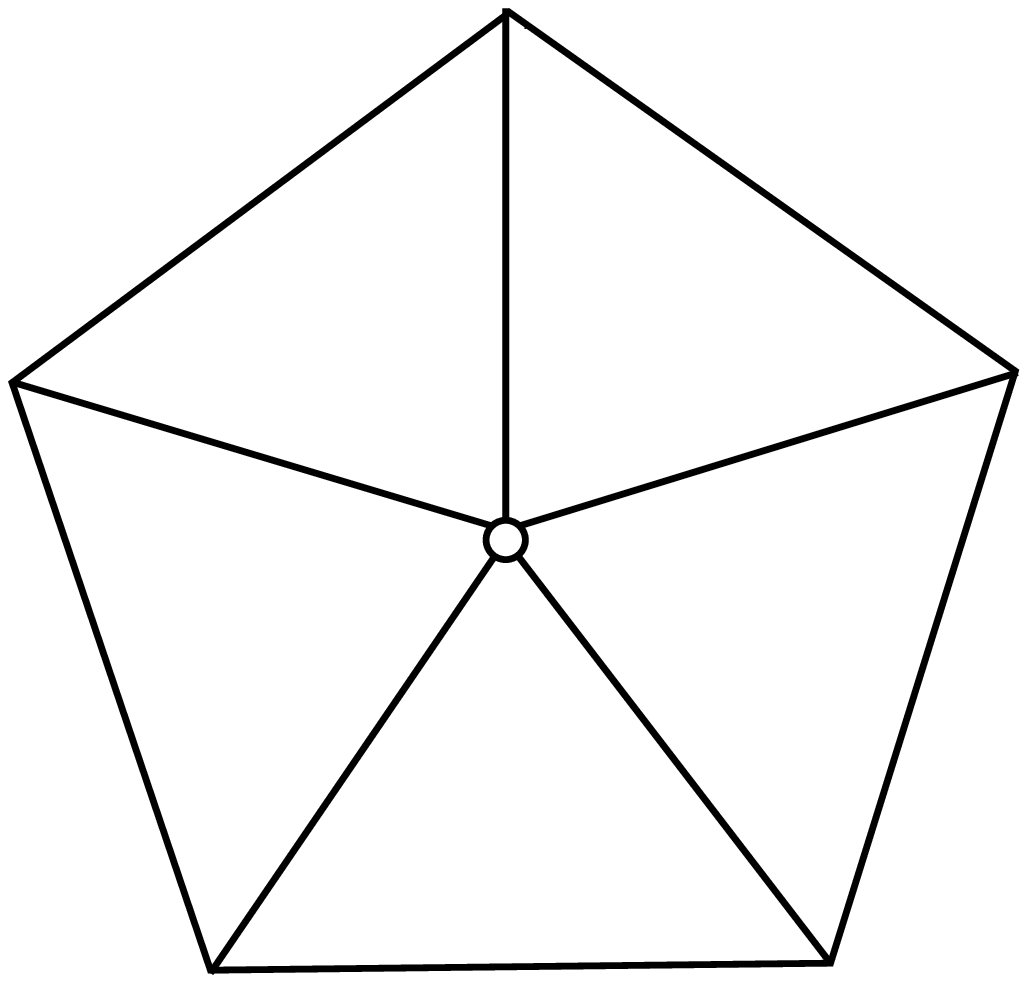}&
5&
1&
$\kappa=\frac{2\pi}{6}=\frac{\pi}{3}$&
$\frac{5\pi}{3}$\\
\hline
\includegraphics[scale=.1]{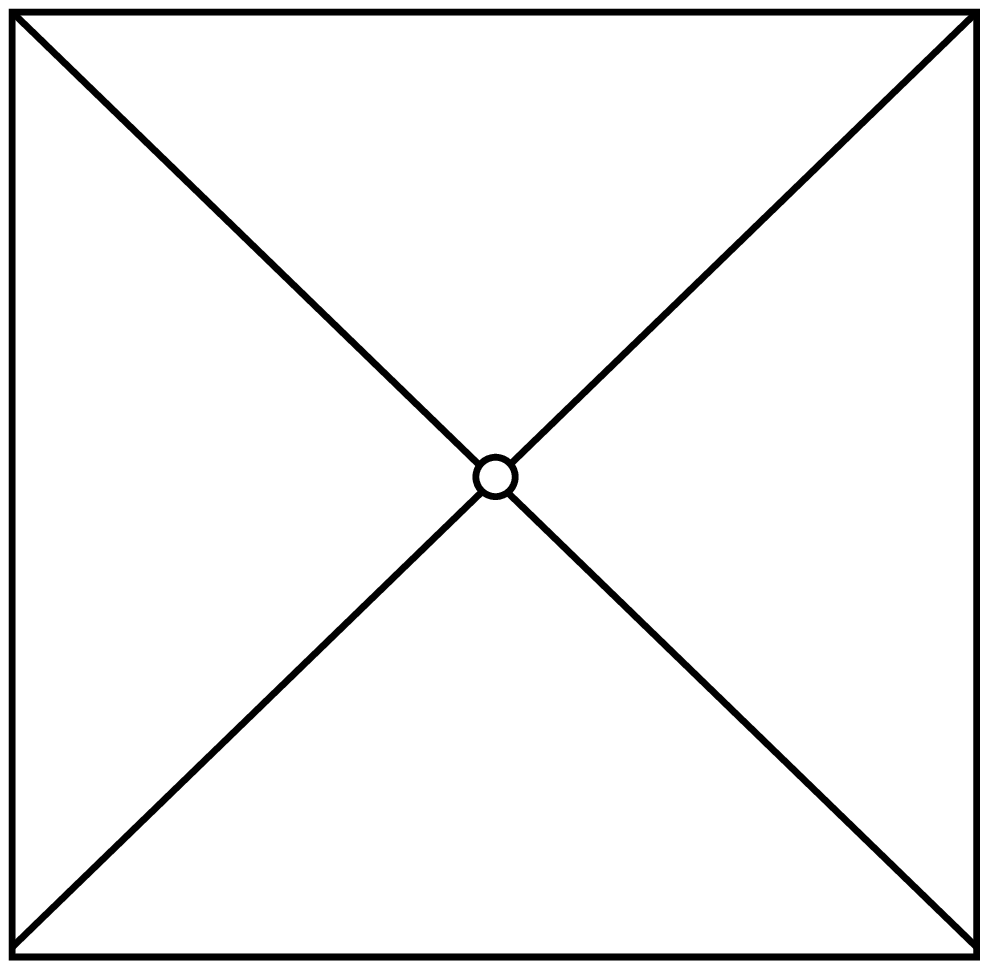}&
4&
2&
$\kappa=\frac{4\pi}{6}=\frac{2\pi}{3}$&
$\frac{4\pi}{3}$\\
\hline 
\includegraphics[scale=.1]{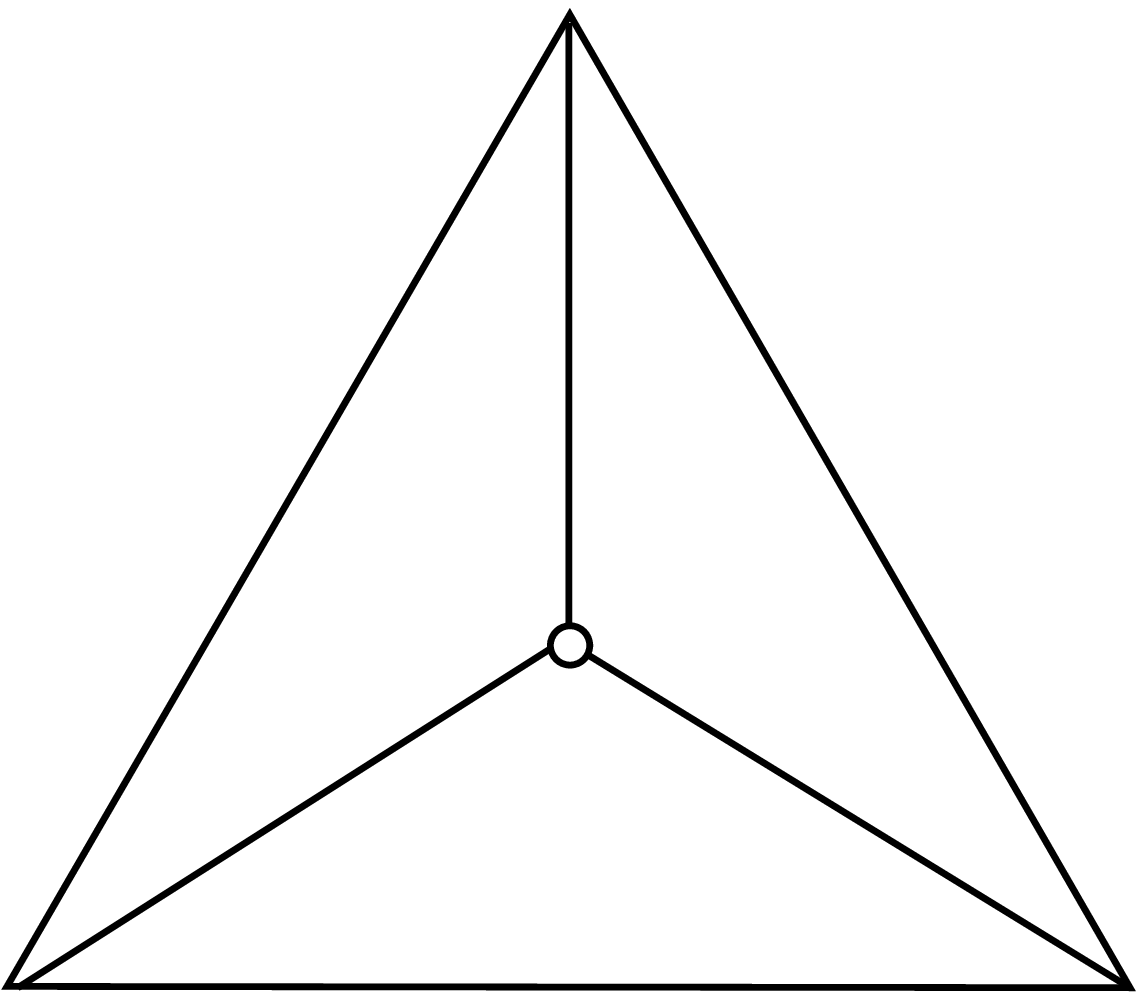}&
3&
3&
$\kappa=\frac{6\pi}{6}=\pi$&
$\pi$\\ \hline 
\includegraphics[scale=.1]{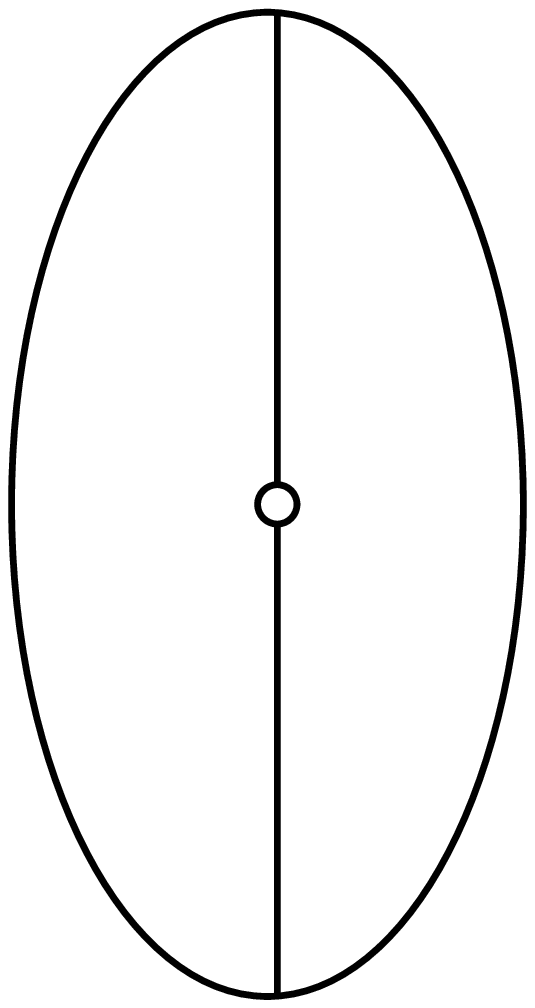}&
2&
4&
$\kappa=\frac{8\pi}{6}=\frac{4\pi}{3}$&
$\frac{2\pi}{3}$\\ \hline 
\includegraphics[scale=.1]{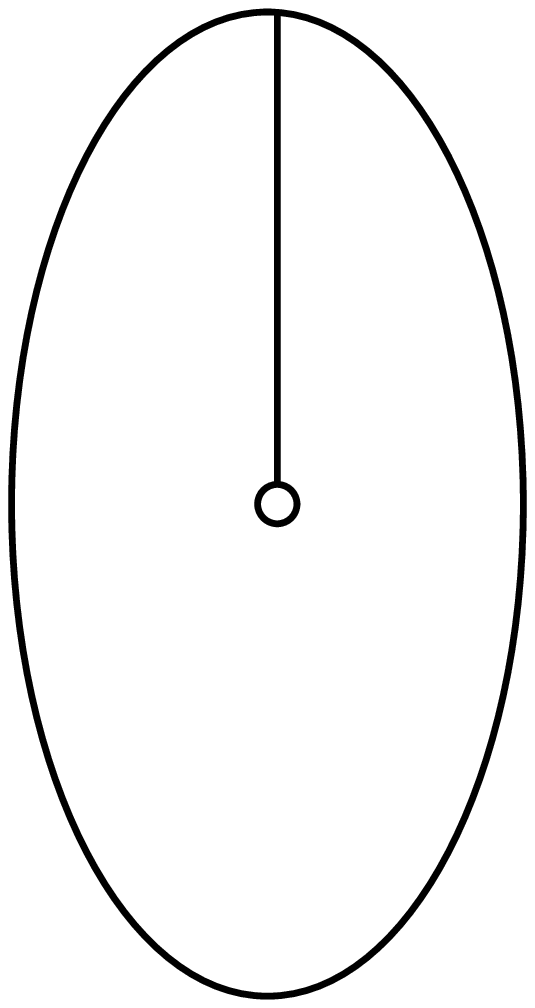}&
1&
5&
$\kappa=\frac{10\pi}{6}=\frac{5\pi}{3}$&
$\frac{\pi}{3}$\\ \hline 
\end{tabular}

\phantom{.}
\bigskip
A vertex of  an equilateral triangulation is said to be {\it non-negatively curved} \index{non-negatively curved vertex}\index{vertex ! non-negatively curved} if
there are at most six triangles meeting at that vertex and {\it positively curved} \index{positively curved vertex}\index{vertex ! positively curved} if there are at
most 5 triangles meeting at that vertex. A triangulation is said to be {\it non-negatively curved} \index{non-negatively curved triangulation} \index{triangulation ! non-negatively curved}
if all its vertices are non-negatively curved. Non-negatively curved triangulations form a very
special class. Basic application of Euler's formula shows that a sphere triangulation of non-negative
curvature may have at most 12 vertices of positive curvature.

Note that one may simultaneously subdivide any Euclidean triangulation, the vertices added in the process will be of zero curvature. Hence we may view these subdivisions as integers rescalings of the original Euclidian structure.

\subsection{Hypergeometric triangulations.}
Thurston studied in the eighties  non-degenerate sphere triangulations of non-negative curvature.
He gave a very concrete and explicit classification and a construction of these sphere
triangulations. These triangulations are related to the works of Picard, Terada, Deligne and Mostow
(PTMD) on higher dimensional hypergeometric functions. It seems appropriate to call these
triangulations {\it hypergeometric} \index{hypergeometric triangulation} \index{triangulation ! hypergeometric}. To any sphere triangulation, there correspond a genus-0
covering of the modular orbifold \index{hypergeometric covering} \index{covering ! hypergeometric} , a modular graph, and a subgroup of the
modular group \index{hypergeometric subgroup} \index{subgroup ! hypergeometric}, each of which we shall call {\it hypergeometric} if the triangulation is
hypergeometric. Recall that the non-degeneracy of the triangulation translates as the absence of terminal edges in the modular graph; or absence of torsion elements in the subgroup and  in this case the covering orbifold is actually a surface.

Thurston showed that hypergeometric triangulations come in (essentially) finitely
many infinite families. These families are parametrized by a finite number of vectors
$\mu=(\mu_1,\cdots, \mu_k)\in \Q_{>0}^k$. The family corresponding to the parameter $(\frac{1}{6},
\frac{1}{6}, \dots)$ of length 12 is the largest family, and all other families can be obtained
from this family by certain degeneration operations. We shall denote by $HG(\mu)$ the family of
hypergeometric coverings related to the parameter $\mu$. One has thus
$$
HG(\mu)\subseteq {\rm {\bf FCov}}_0 (\mathcal M),
$$
where the right-hand side means the genus-0 piece of the covering category.
The parameters $\mu$ also appear in PTMD theory and corresponds to some discrete complex hyperbolic
groups of finite covolume. There is an alternative way to understand these parameters, as we detail in the next section.

There is another way of introducing a flat structure on a curve, via quadrangulations instead of
triangulations. This approach is related to the ${\rm \bf FSub}$ of the group
$\Z\!/\!2\Z*\Z\!/\!4\Z$. 
Quadrangulations are related to the ring of Gaussian integers. Triangulations are related to the
ring of Eisenstein integers. Although it is not explicitly stated in Thurston's paper, one of the
lattices (and its degenerations) he discovered classifies hypergeometric square tilings.

Before going to the heart of the matter, we want to point out  two amusing connections.

\subsection{Fullerenes,  quilts and netballs}
The most famous one among the hypergeometric triangulations \index{icosahedral triangulation} \index{triangulation ! icosahedral} is the icosahedral
triangulation, which belongs to the biggest family of triangulations mentioned above. Many
combinatorial objects with nice properties can be naturally related to hypergeometric
triangulations. They appear spontaneously in diverse fields and there is a very rich terminology
surrounding them. Triangulated spheres are sometimes called {\it deltahedra} \index{deltahedra}. Polyhedra with all
vertices of degree 3 are named {\it trivalent polyhedra}\index{trivalent polyhedra} \index{polyhedra ! trivalent}. In organic chemistry, trivalent polyhedra
with only pentagonal or hexagonal vertices are called {\it fullerenes} \index{fullerenes} (alternative names are:
footballene, buckyballs, buckminsterfullerenes). Fullerenes are studied in chemistry in connection
with the discovery of some complex molecules formed by carbon atoms. In the chemistry literature,
there are catalogs of fullerenes \cite{nomenclature}. Any trivalent polyhedron yields an associated
deltahedron (i.e.  a sphere triangulation) via central subdivision, the associated deltahedron of a
fullerene is then a {\it hypergeometric} triangulation lying in the class which also contains the
icosahedral triangulation. The icosahedron itself corresponds to the molecule $C_{60}$. In the
context of chemistry, coverings in  $HG(\mu)$ with the same branch behavior (passport) appear as
isomers. The question of isomer counting of fullerenes is also being studied in the chemistry
literature. The hypergeometric connection relates this problem to counting orbits of points of a
certain lattice, under the action of a group of automorphisms. 

The fullerenes appear in another, even more surprising context.
Quilts \index{quilts} were invented by Norton to study the ``genus-0 phenomenon" related to the monster group
\cite{hsu}. We may understand quilts as  dessins supplied with some extra information. There is
a special class of quilts, named footballs or netballs by Norton, they appear in the study of
monster and its subgroups. In fact, the netball quilts are precisely fullerenes, and fullerenes are
hypergeometric. It seems that the celebrated genus-0 phenomenon have some connection to hypergeometric
triangulations. An independent sign indicating a possible relevance of hypergeometric triangulations and hyperbolic
geometry to the monster is given by the conjectural ``monstrous proposal" \cite{allcock}.

\section{Branched covers of the sphere.} \label{branched} There is a well-known
classification of branched Galois coverings \index{branched Galois covering} \index{covering ! branched Galois} $\p^1\rightarrow \p^1$; their signature belong to the
list $(m,m)$, $(2,2,m)$, $(2,3,3)$, $(2,3,4)$, $(2,3,5)$. It is also known that the signatures
$(2,3,6)$, $(2,4,4)$ or $(3,3,3)$ and $(2,2,\infty)$ are realized by branched Galois coverings of
$\p^1$ by elliptic curves (or by ${\mathbb A}^1$).

The problem of existence, enumeration and classification of all branched coverings (Galois or not)
of $\p^1$ is an important problem and with the discovery of connection with moduli spaces, 
considerable current research is being devoted to this topic. We may call this bundle of problems 
``the Hurwitz program" \index{Hurwitz program}. This program is of course intractable
in this generality and it is necessary to impose some restrictions, i.e. on the branching behavior
of the coverings. Let us consider the following special instance of the Hurwitz program: 

\medskip\noindent
{\bf Problem E.} Classify all covers
$f:\p^1\rightarrow \p^1$ such that $f$ has ramification index 2 at each fiber above $0\in \p^1$,
ramification index 3 at each fiber above $1\in \p^1$ and has $k_i\geq 0$ points of ramification
index $i$ above $\infty\in\p^1$ for $i=1,2,3,\dots$.

\medskip
We shall see that this problem admits a complete and beautiful solution (by Thurston), 
under the assumption that $k_i=0$ for $i\geq 7$. Obviously, solving it amounts to the 
classification of subgroups of the modular group satisfying a certain regularity condition (of
being genus-0 and torsion-free; equivalently the covering must factor through the covering of the
modular orbifold by $\p^1\backslash\{0,1,\infty\}$). Suppose $f$ is of degree $d$.The Riemann-Hurwitz 
formula \index{Riemann-Hurwitz formula} yields
\begin{equation}\label{orbeuler}
2= e(\p^1)=d\cdot
e(\p^1\backslash\{0,1,\infty\})+\frac{d}{2}+\frac{d}{3}+\sum_{i=1}^\infty
k_i=-\frac{d}{6}+\sum_{i=1}^\infty k_i
\end{equation}
where $e(\p^1\backslash\{0,1,\infty\})=-1$ is the Euler characteristic. Since $\sum_{i=1}^\infty ik_i=d$, one has
\begin{equation}\label{rh}
\sum_{i=1}^\infty (6-i)k_i=12
\end{equation}
The above-mentioned regularity conditions says in effect: the standart triangulation of $\p^1$ with
two triangles having vertices at 0, 1 and $\infty$ lifts to a triangulation of $\p^1$ in a nice
manner. Assume now that $k_i=0$ for $i>6$ and note that the number $k_6$ does not have any effect
in the above formula. According to the terminology of Thurston, the condition $k_i=0$ for $i>6$,
means that the lifted triangulation is of {\it non-negative combinatorial curvature.} \index{non-negative combinatorial curvature} \index{curvature ! non-negative combinatorial}
Quilts satisfying this condition are called {\it 6-transposition quilts}, \index{6-transposition quilts} \index{quilts ! 6-transposition}since the icosahedral quilt is a
football, Norton also suggested the name {\it netballs}\index{netballs} (see \cite{hsu}). We shall simply call
them (be it quilt, triangulation, subgroup or covering): {\it hypergeometric}.

By $[n]_k$ we shall denote a sequence which consists of $k$ repetitions of $n$.

We may present the solutions of (\ref{rh}) subject to the restriction $k_i=0$ for $i>6$ by vectors
$\mu=([1]_{k_1},[2]_{k_2},[3]_{k_3},[4]_{k_4},[5]_{k_5})$ (if we ignore $k_6$ then the list is finite). Let
us denote by $HG_{Eis}(\mu)$ (read as: ``the class of hypergeometric curves of type $\mu$") the
corresponding set of branched coverings, so one has a natural inclusion
$$
HG_{Eis}(\mu)\subseteq {\rm {\bf FCov}}_0 \mathcal M.
$$
A solution of (\ref{rh}) is $([5]_{12})$. If $k_6=0$ it is known that there exists indeed a covering
with this branch data, namely the icosahedral covering of signature $(2,3,5)$. Simultaneous subdivisions are also of the same type. Hence, the set
$HG_{Eis}(\mu)$ is infinite for $\mu=([5]_{12})$. In fact, the set $HG_{Eis}(\mu)$ contains many other elements as we shall see below. 
What is  surprising is that the full set of
solutions of (\ref{rh}) yields exactly those entries in Picard-Terada-Deligne-Mostow's list of
reflection groups that corresponds to Eisenstein integers; these solutions are tabulated in the appendix. 
Apparently, there is an alternative way
of understanding Deligne-Mostow's integrality conditions, which may explain some surprising coincidences appearing in this field.
{\it Notice the change of view here:} 12 moving points of Deligne and Mostow are rigidified and become fibers
above infinity of a covering of the modular orbifold, in other words, cusps of a modular curve.

\subsection{The Gaussian Case.}
The solution of Problem E is connected to Eisenstein integers. \index{Eisenstein integers} There is another problem which admits a similar solution, which is connected to Gaussian integers. \index{Gaussian integers}

\medskip
\noindent
{\bf Problem G.} Classify all covers
$f:\p^1\rightarrow \p^1$ such that $f$ has ramification index 2 at each fiber above $0\in \p^1$,
ramification index 4 at each fiber above $1\in \p^1$ and has $k_i\geq 0$ points of ramification
index $i$ above $\infty\in\p^1$ for $i=1,2,3,\dots$.

\medskip
We shall see that this problem admits a complete and beautiful solution, 
under the assumption that $k_i=0$ for $i\geq 5$.
Obviously, solving it amounts to the 
classification of subgroups of the triangle group $\ztwo\star\!\zfour$,
satisfying a certain regularity condition (of
being genus-0 and torsion-free; equivalently the covering must factor through the (non-Galois) covering of the triangle orbifold of signature $(2,4,\infty)$ by $\p^1\backslash\{0,1,\infty\}$. 
Note that the existence of this covering shows that this triangle orbifold is 
commensurable with the modular orbifold.)
Suppose $f$ is of degree $d$.  The Riemann-Hurwitz formula yields
\begin{equation}\label{orbeuler2}
2= e(\p^1)=d\cdot
e(\p^1\backslash\{0,1,\infty\})+\frac{d}{2}+\frac{d}{4}+\sum_{i=1}^\infty
k_i=-\frac{d}{4}+\sum_{i=1}^\infty k_i
\end{equation}
where $e(\p^1\backslash\{0,1,\infty\})=-1$ is the Euler characteristic.

Since $\sum_{i=1}^\infty ik_i=d$, one has
\begin{equation}\label{rh2}
\sum_{i=1}^\infty (4-i)k_i=8
\end{equation}
The maximal abelian covering of the triangle orbifold of signature $(2,4,\infty)$ is a punctured torus. 
Coverings of the latter orbifold yields quadrangulated surfaces, (or origamis \index{origamis}) which is studied in the context of billards \index{billiards}and in Teichm\"uller theory.
Assume now that $k_i=0$ for $i>5$ and note that the number $k_4$ does not have any effect
in the above formula. The condition $k_i=0$ for $i>6$,
means that the lifted square tiling is of {\it non-negative combinatorial curvature}\index{non-negative combinatorial tiling} \index{tiling ! non-negative combinatorial}.
We shall call these tilings {\it hypergeometric}.
(The class of quadrangulations studied in billards usually possess singularities of negative combinatorial curvature, so they are not hypergeometric in this sense). 

We may present the solutions of (\ref{rh}) subject to the restriction $k_i=0$ for $i>4$ by vectors
$\mu=([1]_{k_1},[2]_{k_2},[3]_{k_3})$ (if we ignore $k_4$ then the list is finite). Let
us denote by $HG_{Gauss}(\mu)$ (read as: ``the class of hypergeometric quadrangulations of type $\mu$") the
corresponding set of branched coverings, so one has a natural inclusion
$$
HG_{Gauss}(\mu)\subseteq {\rm {\bf FSub}}_0^*\, \ztwo * \zfour
$$
where on the right we have the conjugacy classes of finite-index subgroups inside $\ztwo * \zfour$.
A solution of (\ref{rh2}) is $([3]_{8})$. If $k_4=0$ it is known that there exists indeed a covering
with this branch data, namely the tetrahedral covering of signature $(2,4,3)$. Hence, the set
$HG_{Gauss}(\mu)$ is non-empty for $\mu=([5]_{12})$. What is  surprising is that the full set of
solutions of (\ref{rh2}) yields exactly those entries in Picard-Terada-Deligne-Mostow's list of
reflection groups that corresponds to Gaussian integers; these solutions are tabulated below. 

\bigskip
\begin{center}
{\small
\begin{tabular}{|c|ccc|c|c|c|c|c|c|c|}
\hline
dim & $k_1$ & $k_2$ & $k_3$ &\!deg\!&\!\!Compct?\!\!&\!Number\!&\!Pure?\!&\!ar?\!\\\hline
5&0&0&8&2&N&3&P&AR\\\hline
4&0&1&6&2&N&4&P&AR\\\hline
3&1&0&5&2&N&5&P&AR\\\hline
3&0&2&4&2&N&6&P&AR\\\hline
2&1&1&3&2&N&7&P&AR\\\hline
2&0&3&2&2&N&8&P&AR\\\hline
1&2&0&2&-&N&&&AR\\\hline
1&1&2&1&-&N&&&AR\\\hline
1&0&4&0&-&-&self&&AR\\\hline
0&2&1&0&-&-&self&&AR\\\hline
\end{tabular}}
\end{center}

\subsection{Some Numerology.}
The fact that a Hurwitz-type classification problems A and B admits a very nice solution is encouraging.
Can one relax the above-mentioned conditions of regularity to obtain classifications of some new families of triangulations and discover new discrete complex hyperbolic groups  generated by reflections? Let us relax Problem E as follows:

\medskip\noindent
{\bf Problem E$^\prime$.} 
Classify all covers
$f:\p^1\rightarrow \p^1$ such that $f$ has ramification index 2 or 1 at each fiber above $0\in \p^1$,
ramification index 3 or 1 at each fiber above $1\in \p^1$ and has $k_i\geq 0$ points of ramification
index $i$ above $\infty\in\p^1$ for $i=1,2,3,\dots$. 
\medskip

Suppose $f$ is of degree $d$. 
Let $m_i$ be the number of points above $0$ of ramification index $i$ for $i=1,2$. 
Similarly, let $n_i$ be the number of points above $1$ of ramification index $i$ for $i=1,3$. 
Thus, $m_1+2m_2=n_1+3n_3=\sum_{i=1}^\infty ik_i=d$. The Riemann-Hurwitz formula yields
\begin{equation}
2=-d+(m_1+m_2)+(n_1+n_3)+\sum_{i=1}^\infty
k_i=-d+\frac{d+m_1}{2}+\frac{d+2n_1}{3}+\sum_{i=1}^\infty k_i
\end{equation}
Therefore 
$$
2=-\frac{d}{6}+\frac{m_1}{2}+\frac{2n_1}{3}+\sum_{i=1}^\infty k_i,
$$
and setting $d=\sum_{i=1}^\infty ik_i$ yields
\begin{equation}\label{rh3}
\sum_{i=1}^\infty (6-i)k_i=12-3m_1-4n_1.
\end{equation}
The case $m_1=n_1=0$ was considered in Problem E. Assuming that at least one of $m_1$ and $n_1$ is non-zero, we get the table in Appendix 2.

Of special interest are those cases where the number of fibers above is at least five. 
There are 22 of them; they will conjecturally classify some degenerate triangulations  
and yield some lattices. Equivalently, this will give a classification of a certain family of 
subgroups in the modular group, of genus 0 and with some torsion. 
There is a possibility that these lattices are all commensurable with those in the PTDM list.

\begin{figure}\label{devils}
\begin{center}
\resizebox{0.8\hsize}{!}{\includegraphics{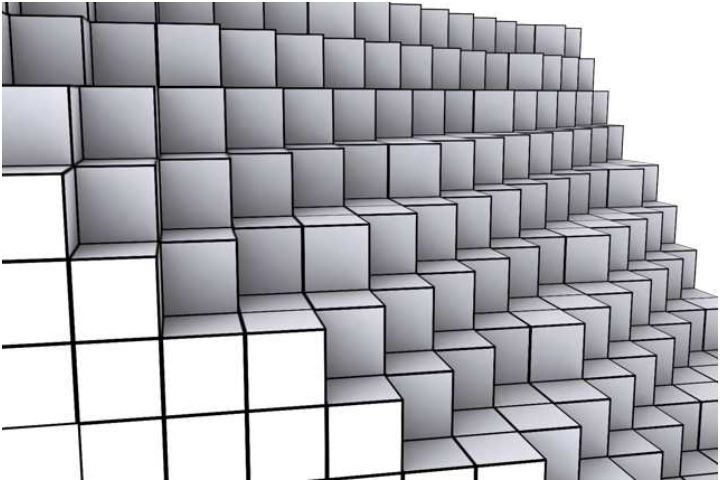}}
\caption{A non-hypergeometric quadrangulation - which points are of positive, zero and negative curvature?  {\tiny{(graphics: courtesy of Mick West)}}}
\end{center}
\end{figure}

In the Gaussian case, one has an analogous modification. 

\medskip\noindent
{\bf Problem G$^\prime$.} 
Classify all covers
$f:\p^1\rightarrow \p^1$ such that $f$ has ramification index 2 or 1 at each fiber above $0\in \p^1$,
ramification index 4, 2 or 1 at each fiber above $1\in \p^1$ and has $k_i\geq 0$ points of ramification
index $i$ above $\infty\in\p^1$ for $i=1,2,3,\dots$. 
\medskip

Suppose $f$ is of degree $d$. 
Let $m_i$ be the number of points above $0$ of ramification index $i$ for $i=1,2$. 
Similarly, let $n_i$ be the number of points above $1$ of ramification index $i$ for $i\in\{1,2,4\}$. 
Thus, $m_1+2m_2=n_1+2n_2+4n_4=\sum_{i=1}^\infty ik_i=d$.The Riemann-Hurwitz formula yields
\begin{equation}
2=-d+(m_1+m_2)+(n_1+n_2+n_4)+\sum_{i=1}^\infty
k_i.
\end{equation}
Therefore 
$$
2=-\frac{d}{4}+\frac{m_1}{2}+\frac{3n_1+2n_2}{4}+\sum_{i=1}^\infty k_i,
$$
and setting $d=\sum_{i=1}^\infty ik_i$ yields
\begin{equation}
\sum_{i=1}^\infty (4-i)k_i=8-2m_1-3n_1-2n_2.
\end{equation}
The case $m_1=n_1=n_2=0$ was considered in Problem G. 
Assuming that at least one of $m_1$, $n_1$ and $n_2$ is non-zero, we get the table in Appendix 3.

\begin{figure}[h]
\begin{center}
\resizebox{0.4\hsize}{!}{\includegraphics{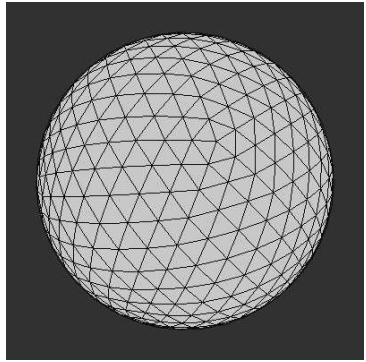}}
\caption{A hypergeometric sphere triangulation}
\end{center}
\end{figure}

\section{Thurston's work on sphere triangulations} \label{thurston} 
We must stress that the lists of the previous section are purely hypothetical. Numerology exhibits potentialities but doesn't say anything about their realizations. Attacking this problem in a straightforward manner requires studying monodromy presentations, 
which is a time and space consuming combinatorial problem that one may hope to attack by a computer.
In contrast with this, one of the results stated in Thurston's 1987 preprint is the following theorem

\medskip\noindent
{\bf Theorem (Thurston, \cite{shapes})} {\it (Polyhedra are lattice points) There is a lattice
${\cal L}$ in complex Lorenz space $\CC^{(1,9)}$ and a group $\Gamma_{DM}$ of automorphisms, such
that sphere triangulations of non-negative combinatorial curvature are elements of ${\cal L}_+/\Gamma_{DM}$, where
${\cal L}_+$ is the set of lattice points of positive square-norm.
The square norm of a lattice
point is the number of triangles in the triangulation. The projective action of $\Gamma_{DM}$ on complex projective hyperbolic space $\CC{\mathbb H}^9$ (the unit ball in
$\CC^9\subset \CC{\mathbb \p}^9$) has quotient of finite volume. }

\medskip\noindent
This lattice was explicitly identified by Allcock \cite{allcock}. 
Triangulations lying on the same line through the origin are simultaneous subdivisions of a ``{primitive}" triangulation on the line and therefore define isometric polyhedra. Hence the projectification  
$$
\mathbb P{\cal L}_+/\Gamma_{DM}\subset {\cal M}_{DM}:=\CC{\mathbb H}^9/\Gamma_{DM}
$$ 
classifies the isometry classes (``{shapes}" in Thurston's terms) of polyhedra, where ${\cal M}_{DM}$ is the ball-quotient space $\CC{\mathbb H}^9/\Gamma_{DM}$. We shall call these ``{hypergeometric points}" of the moduli space.
Thurston also describes a very explicit method to construct these triangulations and gives the estimation $O(n^{10})$
for the number of triangulations in $HG_{Eis}([1]_{12})$ with up to $2n$ triangles. 

\medskip\noindent
{\bf Problem:} (Isomer counting) 
Let $\Delta_n(\mu)$ be the number of triangulations in $HG(\mu)$ with $n$ triangles. 
Find an appropriate generating function $T_\mu$ for the numbers $\Delta_n(\mu)$.

\medskip
It must be possible to complete Thurston's results as follows:

\medskip\noindent
{\bf Theorem.} {\it Let $\mu$ be an admissible curvature vector of length $\ell(\mu)=\ell$. There is a lattice
${\cal L(\mu)}$ in complex Lorenz space $\CC^{(1,\ell-3)}$ and a group $\Gamma(\mu)$ of automorphisms, such
that triangulations of type $\mu$ are elements of ${\cal L}_+(\mu)/\Gamma(\mu)$, where
${\cal L}_+(\mu)$ is the set of lattice points of positive square-norm. The projective action of
$\Gamma{(\mu)}$ on complex projective hyperbolic space $\CC{\mathbb H}^{\ell-3}$ (the unit ball in
$\CC^{\ell-3}\subset \CC{\mathbb \p}^{\ell-3}$) has quotient of finite volume. The square norm of a lattice point is the number of triangles in the triangulation.}

\medskip
The previous theorem corresponds to the longest parameter $\mu=[1]_{12}$, and the other 
pairs $({\cal L(\mu)},\Gamma(\mu))$ arise as degenerations of this one. As abstract groups, 
$\Gamma{(\mu)}$ are braid group quotients. We denote the quotient 
$$
\CC{\mathbb H}^{\ell-3}/\Gamma(\mu)=:\mathcal M_\mu.
$$ 
As above, there is a dense subset of hypergeometric points inside the ball quotient space $M_\mu$:
$$
\mathbb P{\cal L_+(\mu)}/\Gamma(\mu)\subset {\cal M}_{\mu}.
$$
These points are conjecturally defined over $\overline{\mathbb Q}$.
It is an important task to understand the structure of the ``hypergeometric web" \index{hypergeometric web} \index{web ! hypergeometric}, i.e. various degenerations of triangulations in this 9-dimensional moduli space (with respect to the  Galois action). Even the integral lattices themselves have not been explicitly identified in the literature. \.Ismail Sa\u glam \cite{ismailthesis}, \cite{ismail} proved this theorem for the cases 
$\mu=[2]_{6}$ and $\mu=[3]_{4}$ (and also $\mu=[1]_{12}$),
using alternative and more explicit methods than Thurston's hard-going paper.
His proof gives a construction of those triangulations and also applies to Ayberk Zeytin's theorem concerning quadrangulations presented below. In case  $\mu=[3]_{4}$, the group in question is the modular group, i.e.
 $\Gamma([3]_{4})\simeq \mathrm{PSL}_2(\Z)$ and provides the most amenable family of triangulations and polyhedra on which the Galois action should be studied. We shall give a construction of this family in the last section of the current paper.

Allcock gave in the late 1990's a more direct construction of $\Gamma_{DM}$ as a group of automorphisms
of the lattice ${\cal L}$ and imitated this construction to build a 13-dimensional ball quotient
related to a lattice ${\cal L}^A$ which is derived from the Leech lattice \cite{allcock}. His construction is
conjecturally related to the Monster group in a precise way \cite{monstrous}. The connection we
unearthed above between the hypergeometric triangulations and the quilts related to the monster
(see \cite{hsu}, Chapter 11) reveals that there is something about the monster in the
hypergeometric world. Is there a similar combinatorial interpretation of Allcock's lattice ${\cal
L}^A$ i.e. as a set of triangulations? If yes, most of the questions we raise here about the
Deligne-Mostow's ball quotients and related objects could be formulated for Allcock's ball-quotient
as well.

As for the quadrangulations, one has the following result

\medskip\noindent
{\bf Theorem} (Ayberk Zeytin \cite{ayberkthesis}, \cite{ayberkmakale}) {\it (Quadrangulations are lattice points) There is a lattice
${\cal L}$ in complex Lorenz space $\mathbb C^{(1,8)}$ and a group $\Gamma_{DM}$ of automorphisms, such
that quadrangulations of non-negative combinatorial curvature are elements of ${\cal L}_+/\Gamma_{DM}$, where
${\cal L}_+$ is the set of lattice points of positive square-norm. The projective action of
$\Gamma_{DM}$ on complex projective hyperbolic space $\CC{\mathbb H}^9$ (the unit ball in
$\CC^9\subset \CC{\mathbb \p}^9$) has quotient of finite volume. The square of the norm of a lattice
point is the number of quadrangles in the triangulation.}

\begin{figure}
\begin{center}
\resizebox{0.3\hsize}{!}{\includegraphics{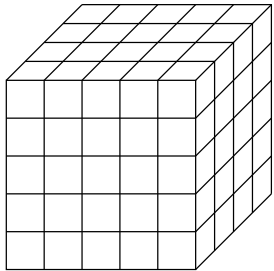}}
\caption{A hypergeometric sphere quadrangulation}
\label{bigcube}
\end{center}
\end{figure}

\subsection{Hypergeometric functions, ball-quotients of Picard, Terada, Deligne and Mostow and
the transcendence results of Wolfart and Shiga.} \label{hyper} 
Multivariable hypergeometric functions \index{multivariable hypergeometric function} \index{hypergeometric function ! multivariable} arise as the uniformization maps of the moduli spaces ${\mathcal M}_\mu$.

The hypergeometric differential equation \index{hypergeometric differential equation} \index{differential equation ! hypergeometric}
was first discovered by Euler, the term {\it hypergeometric} is even older; the name {\it Gauss'
hypergeometric functions} \index{Gauss' hypergeometric functions} \index{hypergeometric function ! Gauss'} is also frequently used after Gauss's contributions. Appell introduced a
two-variable hypergeometric function, which was further generalized to arbitrary many variables by
Lauricella. Following the works of Riemann and Schwarz in dimension one, Picard studied the
finiteness and discreteness of monodromy for Appell's hypergeometric functions. Terada extended this
work to the Lauricella hypergeometric functions \index{Lauricella hypergeometric function} \index{hypergeometric function, Lauricella} in the 1970's. Deligne and Mostow's paper on Lauricella
hypergeometric functions appeared in the1980's and gave a uniform and rigorous treatment of
discreteness using algebraic geometry (see \cite{looi} for an elementary treatment). Thurston used
geometric methods to reprove these discreteness results, without mentioning hypergeometric
functions at all \cite{shapes}. By using the numerical ball-quotient criterion (Miyaoka-Yau
proportionality) Hirzebruch,  Holzapfel  and followers discovered some other
discrete complex hyperbolic groups generated by reflections, but they all turned out to be
commensurable with a lattice in Picard-Terada-Deligne-Mostow's (PTMD) list \cite{commens}.
Recently, Heckman-Couwenberg-Looijenga gave another generalization and obtained some other complex
hyperbolic reflection groups \cite{hcl}. However, it is not known if these lattices are
commensurable with the PTMD lattices. Yoshida and collaborators gave alternative modular
interpretations of hypergeometric functions and studied their properties \cite{mylove}.

Transcendence problems for the (multivalued) Schwarz maps \index{Schwarz maps}have been studied by Cohen, Wohlfart, Shiga and
Suzuki, their result for higher dimensions roughly reads: ``if the Schwarz map value at $\tau$ is
algebraic, then a certain Prym variety parametrized by $\tau$ has CM". On the reverse direction, it
is a natural wonder what the images of lattice points under the ball-quotient maps (inverse
Schwarz maps) are. Is it possible to compute their precise values? We conjecture
that the images of lattice points are dense, and algebraic. Finally, the Galois action 
is compatible with the action on the corresponding
hypergeometric curves. Moreover, the Galois action must respect the structure of the
``hypergeometric web", which is formed by the degenerations in the 9-dimensional ancestral
ball-quotient.
\subsection{Questions.}
In the light of their connections to hypergeometric functions, combinatorics and group theory,
there is a well-founded hope that hypergeometric curves will render themselves to explicit calculation
and unlike the case of a general dessin, we can study the Galois action on them.
{There are several circles of questions that appear:} 

By ``hypergeometric triangulation or quadrangulation" \index{hypergeometric quadrangulation} \index{quadrangulation ! hypergeometric} (equivalently ``hypergeometric dessin") of type $\mu$, we both mean a point in ${\mathcal L}(\mu)$ and the sphere triangulation defined by this point. 
``Hypergeometric curve"  \index{hypergeometric curve} \index{curve ! hypergeometric}(or ``hypergeometric cover") of type $\mu$ means the covering of the Riemann sphere defined by a hypergeometric triangulation of type $\mu$. A ``{hypergeometric point}" \index{hypergeometric point} \index{point ! hypergeometric} of type $\mu$ is an element of 
$\mathbb P{\cal L_+(\mu)}/\Gamma(\mu)\subset {\cal M}_{\mu}$, in other words it is a shape parameter of a polyhedra.
Every hypergeometric point represents a ray of hypergeometric triangulations, all obtained from a basic triangulation by simultaneous subdivision.
\subsubsection{Group theory and combinatorics.}
Which hypergeometric curves are modular (i.e. dominated by congruence modular curves)? Given a
hypergeometric curve, find the smallest Galois cover that dominates it. Characterize the monodromy
groups of hypergeometric covers. Compare these monodromy groups with nilpotent and solvable groups;
are these groups non-abelian in an essential manner? Given two hypergeometric covers, find the
(dessin of) smallest covering  that dominates both. Find also the smallest Galois covering that
dominates both. Find an appropriate generating function for the number of hypergeometric triangulations of the same
type $\mu$ (isomer counting).
\subsubsection{Field theory and Galois action.}
Given a hypergeometric triangulation, describe the corresponding Belyi map explicitly and study the
Galois action. Are the Galois action on ${\cal L_+(\mu)}/\Gamma(\mu)$ (defined via hypergeometric curves)  
and the Galois action on the hypergeometric points $\mathbb P{\cal L_+(\mu)}/\Gamma(\mu)$ compatible?
Does this action respect degeneration of triangulations? 
Is the Galois action faithful on hypergeometric curves? (probably it isn't).

Describe the fields of definitions of hypergeometric covers. Describe the minimal field of
definitions $F_\mu^d$ of hypergeometric covers of the same type $\mu$ and degree $d$, and estimate
the order of growth of $[F_\mu^d:\Q]$ as $d\rightarrow \infty$. Describe the minimal field of
definition of all hypergeometric covers of the same type $\mu$. Characterize the minimal field of
definition of all hypergeometric covers.  
\subsubsection{Moduli space, transcendence, rational point counting.}
Show that the hypergeometric points are dense.
Calculate some hypergeometric points explicitly. Is it possible to obtain a triangulation represented by a hypergeometric point ?
Describe the fields of definitions of hypergeometric points. 
Give examples of non-hypergeometric algebraic points of ${\cal M}_{\mu}$.

The minimal number of triangles of hypergeometric triangulations represented by a hypergeometric
point $p_{\cal T}$ defines a ``height" function on the points $p_{\cal T}$. Describe the minimal
field of definitions $K_\mu^d$ of hypergeometric points of the same type $\mu$ and height, and
estimate the order of growth of $[K_\mu^d:\Q]$ as $d\rightarrow \infty$. Describe the minimal field
of definition of all hypergeometric points of the same type $\mu$. Characterize the minimal field of
definition of all hypergeometric points. Count the hypergeometric points. 
\subsubsection{Moonshine. \index{moonshine}}
Elucidate the connections between the netballs of Norton  (group theory), triangulations of
non-negative curvature of Thurston (geometry), hypergeometric curves (algebraic geometry) and
 Allcock's ``monstrous proposal" (complex hyperbolic geometry).
We invite you to inspect the quilts in \cite{hsu} to realize that they are all hypergeometric.
\subsubsection{Hypergeometric  Grothendieck-Teichm\"uller Theory.\index{Grothendieck-Teichmuller Theory}}%
The $\overline{{\cal M}}_{DM}$ includes all families of hypergeometric triangulations as
degenerations. Let us call this structure the ``hypergeometric web". Devise a hypergeometric
version of the Grothendieck-Teichm\"uller group $\widehat{GT}$, deduced from the relations of the
``hypergeometric web" instead of the greater Teichm\"uller tower. 
\subsubsection{Other lattices.}
Thurston's article includes more lattices than those classifying the hypergeometric triangulations
and hypergeometric square tilings. Is there a combinatorial interpretation of these lattices,
similar to triangulations or tilings? Are these lattices connected to some arithmetic curves in
some other way? Do they admit a Galois action?

\subsection{Hypergeometric completion of the profinite modular group.} \label{completion}
Let $G$ be a finitely presented group and let $\widehat{G}$ be its profinite completion. Let
${\mathbb H}=\{H_\alpha\}_{\alpha\in I}$ be a system of finite index subgroups of $G$, satisfying
the property:

\smallskip \noindent
(*) for any $i\in \Z_{>0}$, there are only a finite number of $\alpha$'s such that $[G:H_\alpha]\leq i$.

\smallskip\noindent
To ${\mathbb H}$, one may of associate a quotient $\widehat{G}_{\mathbb H}$ of $\widehat{G}$ as follows:
Let $\overline{\mathbb H}:=\{\overline{H}_\alpha\}_{\alpha\in I}$, where
$\overline{H}_\alpha:=\bigcap_{g\in G}gH_\alpha g^{-1}$ is the normal core of $H_\alpha$.
Then $\overline{\mathbb H}$ also satisfies the property (*), and the normal subgroups
$$
H(i):=\bigcap_{[G:H_\alpha]\leq i} H_\alpha
$$
are of finite index in $G$ as well. Then $G \rhd H(1) \rhd H(2) \rhd \dots $
is a chain of normal subgroups of $G$. Put
$\widehat{G}_{\mathbb H}:=\lim_{\leftarrow} G/H(i)$.
One may call $\widehat{G}_{\mathbb H}$: ``the completion of $G$ with respect to the system ${\mathcal
H}$''. Any system  ${\mathbb H}$ can be enriched by the set of all co-nilpotent (or co-pro-$\ell$,
or co-solvabe.) normal subgroups of all elements in ${\mathbb H}$, yielding a greater system and
an induced ``enriched" completion.

If we take $G=\mathrm{PSL}(2,\Z)$ and ${\cal H}=HG(\mu)$, then the above procedures yield completions
(``enriched" if we wish) $\widehat{\mathrm{PSL}}^\mu(2,\Z)$. This is a somewhat artificial construction, but
it seems that this is the only algebraic object at our immediate disposal, which is derived from
hypergeometry and on which we can study the Galois representation (and not merely a Galois action).
Questions: Is $\widehat{\mathrm{PSL}}^\mu(2,\Z)$ metamotivic? \index{metamotivic} meaning: is it essentially ``{larger}" from
almost nilpotent completions? What is the kernel of the corresponding Galois representation? Can we
get an analogue of the Grothendieck-Teichm\"uller group by considering the total structure of the
hypergeometric web?

\section{Case study: the simplest families of triangulations and quadrangulations}
Here we give an overview of some results from the second named author's thesis \cite{ismailthesis} to describe the family of triangulations $HG_{Eis}([3]_4)$ and the family of quadrangulations  $HG_{Gauss}([2]_4)$. 
The set $HG_{Eis}([3]_4)$ is the set of triangulations with 
4 singular vertices (vertices of non-zero combinatorial curvature) such that each of these vertices is incident to 3 triangles. Similarly, $HG_{Gauss}([2]_4)$ is the set of quadrangulations with 4 singular vertices such that each singular vertex is incident to 2 quadrangles.

We need to introduce some terminology from the theory of cone metrics on 2-dimensional surfaces. 
\subsection{Cone Metrics on Surfaces}
Our reference in this section is \cite{troyanov} and \cite{troyanov1}. 
A {\it triangulated surface} is roughly a surface with an Euclidean triangulation \index{Euclidean triangulation} on it. Here is the formal definition.
\begin{definition}
A \textit {triangulated surface} \index{triangulated surface} \index{surface ! triangulated} is a surface $S$ together with a set of pairs $\mathbf{T}=\{(T_{\alpha},f_{\alpha})\}_{\alpha \in A}$ where each $T_{\alpha}$ is a compact subset of $S$
and each $f_{\alpha}:T_{\alpha}\rightarrow \mathbb{R}^2$ a  diffeomorphism with a non-degenerate euclidean triangle such that
\begin{itemize}
\item[\ITEM]
$T_{\alpha}$'s cover $S$.
\item[\ITEM]
If $\alpha \neq \beta$ then intersection of $T_{\alpha}$ and $T_{\beta}$ is either empty or edge or a vertex.
\item[\ITEM]
If $T_{\alpha}\cap T_{\beta}$ is not empty then there is an element $g_{\alpha \beta}$  $\in E(2)$ (the group of isometries Euclidean plane) such that
$f_{\alpha}=g_{\alpha \beta}f_{\beta}$.
\end{itemize}
\end{definition}
\begin{definition}
A \textit{cone metric} \index{cone metric} on a  triangulated surface is a metric obtained by using given triangulation.
\end{definition}
A surface with a cone metric will be called {\it flat surface}\index{flat surface}\index{surface ! flat}. It is clear that for each point $p$ on a flat surface $S$ there is a notion of angle, $\theta_p$.
The value $\kappa_p=2\pi - \theta_p$ is called the {\it curvature at $p$}. With this preparation we may present the Gauss-Bonnet Theorem \index{Gauss-Bonnet theorem}  \index{theorem ! Gauss-Bonnet}and Hopf-Rinow Theorem \index{Hopf-Rinow theorem} 
\index{theorem ! Hopf-Rinow} for flat surfaces:
\begin{theorem}(Gauss-Bonnet)
Denote by $\chi(S)$ the Euler charateristic of $S$. For any compact flat surface $S$ without boundary we have
\begin{align*}
\sum_{p \in S}(2\pi-\theta_p)=\chi(S).
\end{align*}
\end{theorem}
This formula is easily established by summing angles at singular vertices and counting number of triangles used.
\begin{theorem}(Hopf-Rinow)
Let $S$ be a complete, connected, flat surface. Then any two points in $S$ can be joined by a shortest geodesic in $S$. 
\end{theorem}
How can we obtain cone metrics on sphere? To be more precise, assume that  we are given positive numbers $\theta_1, \theta_2, \theta_3$ so that 
\begin{equation}
\theta_1+\theta_2+\theta_3= 2\pi.
\end{equation}
Can we find a cone metric with 3 singular points such that cone angles at these points are $\theta_1, \theta_2$ and $\theta_3$? Answer for this question is affirmative,
see Figure ~\ref{3sing}.
\usetikzlibrary{calc}
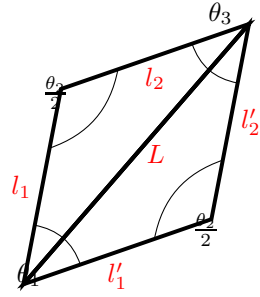
\begin{figure}
\centering

\begin{tikzpicture} 
 \coordinate (A) at  (1.5,{0.5*sqrt(3)})  ;
    \coordinate (B) at (2,{2*sqrt(3)}) ;
    \coordinate (C) at (4,{sqrt(3)});
     \coordinate (D) at (4.5,{2.5*sqrt(3)}) ;
               \
    \draw [black,line width=1.5pt]  (A) -- (B) node[red,midway,left] {$l_1$} -- (B) node[black,midway,left] {$$}-- (D)  node[red,midway,below] {$l_2$}        -- cycle ;   
        \draw [black,line width=1.5pt] (A) -- (C) node[,red,midway,below] {$l_1'$} -- (C) node[black,midway,below] {$$}    -- (D)  node[red,midway,right] {$l_2'$}    -- cycle;  
  \draw [black,line width=1.5pt] (A) --( A)  node[black,midway,below] {$$} -- (D)  node[red,midway,right] {$L$} -- (D)  node[black,midway,right] {$$}  ;     
   
\draw[black!100] (A) -- ($(A)!8mm!(B)$) to[bend left] ($(A)!8mm!(C)$)   node[black,midway,below] {$\theta_1$} -- cycle; 
     
\draw[black!100] (B) -- ($(B)!8mm!(A)$) to[bend right] ($(B)!8mm!(D)$)   node[black,midway,above] {$\frac{\theta_2}{2}$}  -- cycle;                       
\draw[black!100] (D) -- ($(D)!8mm!(B)$) to[bend right] ($(D)!8mm!(C)$) node[black,midway,above] {$\theta_3$}   -- cycle;

\draw[black!100] (C) -- ($(C)!8mm!(A)$) to[bend left] ($(C)!8mm!(D)$)   node[black,midway,right] {$\frac{\theta_2}{2}$}  -- cycle;                    
                        
                         \end{tikzpicture}
    \caption{Constructing a cone metric with 3 singular points} 
    \label{3sing}
    \end{figure}
    In Figure ~\ref{3sing} lengths of $l_1$ and $l_1'$ are equal. Also lengths of $l_2$ and $l_2'$ are equal. If we glue $l_1$ with $l_1'$ and $l_2$
    with $l_2'$ , we get a cone metric on sphere with desired properties. Indeed, this is the only cone metric with above property up to homothety and orientation preserving 
    isometry.
    
    At this point, it is natural to ask whether every cone metric on sphere can be obtained from a polygon in Euclidean plane by identifying some of its edges 
    appropriately. This is not possible in general. However,
     if all curvatures at singular points are positive, the answer is affirmative and is given by Alexandrov 
    Unfolding Process\index{Alexandrov Unfolding Process}.
    \subsubsection{Alexandrov Unfolding Process.}
    Let $\mu$ be a cone metric on sphere with $n$ ($n>2$) singular points of positive curvature. Call these singular points $v_1, v_2, \dots, v_n$. 
Let $s_i$ ($2\leq i \leq n$) be a length minimizing geodesic joining $v_1$ to $v_i$. These geodesics exists by Hopf-Rinow Theorem. It is well known  that $s_i$ and $s_j$ intersect at only $v_1$ when $i \neq j$. If we cut sphere along $s_i$'s, then we can unfold it to the plane without overlapping  as a polygon with $2n-2$ vertices.
Resulting polygon $P$ has $n-1$ vertices coming from $v_1$ and $n-1$ vertices corresponding to $v_i$'s ($ i >1$).  If we glue edeges of this polygon appropriately we get a cone metric on sphere with $n$ singular points. Indeed, this metric, after some normalization, is nothing else than $\mu$. 

This process, Alexandrov Unfolding, briefly says that any cone metric of the positive curvature on sphere can be obtained from a special type of polygon in the plane.
\subsection{Triangulations}
Up to now, we have talked about cone metrics. Now we start to investigate triangulations of sphere. We consider, following \cite{shapes}, a triangulation as a cone metric by assuming that each triangle in triangulation is Euclidean equilateral triangle of edge length 1. We say that two triangulations are equivalent if corresponding metrics are isometric by an orientation preserving isometry sending edges, vertices and triangles to edges, vertices and triangles, respectively. 

How can we construct sphere triangulations? We don't have any means of constructing and classifying them in a systematic manner, other then drawing them by hand. So let us ask a simpler question:  
How can we obtain elements in $HG_{Eis}([4]_3)$\index{$HG_{Eis}([4]_3)$}?
 \usetikzlibrary{calc}
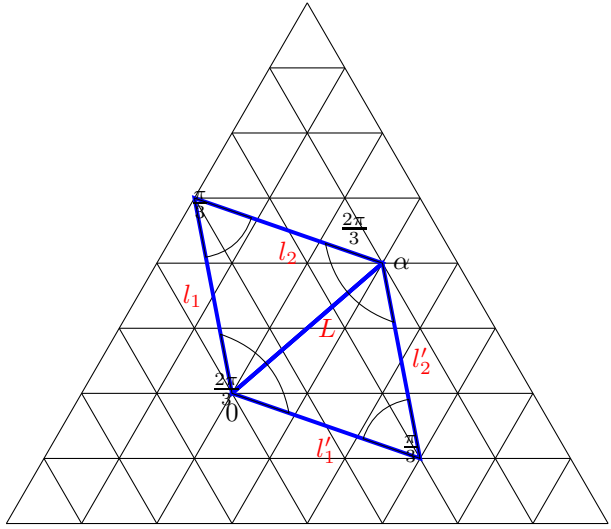
\begin{figure}
\newcommand*\rows{8}
 \centering

\begin{tikzpicture} 
 \coordinate (A) at  (3,{sqrt(3)})  ;
    \coordinate (B) at (2.5,{2.5*sqrt(3)}) ;
    \coordinate (C) at (5.5,{0.5*sqrt(3)});
     \coordinate (D) at (5,{2*sqrt(3)}) ;
               \foreach \row in {0, 1, ...,\rows} {
            \draw ($\row*(0.5, {0.5*sqrt(3)})$) -- ($(\rows,0)+\row*(-0.5, {0.5*sqrt(3)})$) ;    
            \draw ($\row*(1, 0)$) -- ($(\rows/2,{\rows/2*sqrt(3)})+\row*(0.5,{-0.5*sqrt(3)})$);
        \draw ($\row*(1, 0)$) -- ($(0,0)+\row*(0.5,{0.5*sqrt(3)})$);
    }
    \draw [blue,line width=1.5pt]  (A) -- (B) node[red,midway,left] {$l_1$}    -- (D)  node[red,midway,below] {$l_2$}        -- cycle ;   
        \draw [blue,line width=1.5pt] (A) -- (C) node[,red,midway,below] {$l_1'$}    -- (D)  node[red,midway,right] {$l_2'$}    -- cycle;  
  \draw [blue,line width=1.5pt] (A) --( A)  node[black,midway,below] {$0$} -- (D)  node[red,midway,right] {$L$} -- (D)  node[black,midway,right] {$\alpha$}  ;

\draw[black!100] (A) -- ($(A)!8mm!(B)$) to[bend left] ($(A)!8mm!(C)$)   node[black,midway,below] {$\frac{2\pi}{3}$} -- cycle; 
     
\draw[black!100] (B) -- ($(B)!8mm!(A)$) to[bend right] ($(B)!8mm!(D)$)   node[black,midway,above] {$\frac{\pi}{3}$}  -- cycle;                       
\draw[black!100] (D) -- ($(D)!8mm!(B)$) to[bend right] ($(D)!8mm!(C)$) node[black,midway,above] {$\frac{2\pi}{3}$}   -- cycle;

\draw[black!100] (C) -- ($(C)!8mm!(A)$) to[bend left] ($(C)!8mm!(D)$)   node[black,midway,right] {$\frac{\pi}{3}$}  -- cycle;                                  
                        
                         \end{tikzpicture}
    \caption{Obtaining an element of $HG_{Eis}([4]_3)$ from a lattice point.} \label{3triangle}
    \end{figure}
    
Let $\mathbf{Eis}$ \index{Eisenstein integers} be the ring of Eisenstein integers. Observe that $\mathbf{Eis}$ gives a triangulation of the plane. We will obtain desired triangulations from this triangulation. 
Consider the polygon in Figure ~\ref{3triangle} with the following properties:

\begin{itemize}
\item
vertices of the polygon are in $\mathbf{Eis}$,
\item
lengths of $l_1$ and $l_1'$ are equal,
\item
lengths of $l_2$ and $l_2'$ are equal,
\item
angles at $\alpha$ and origin are $\frac{2\pi}{3}$,
\item
angles at the other two vertices are $\frac{\pi}{3}$.
\end{itemize}
If we glue $l_1$ with $l_1'$ and $l_2$ with $l_2'$ we get a triangulation of sphere. Moreover the vertices corresponding to $0$ and $\alpha$ are incident to two triangles.
Also the other two vertices form a single vertex of the triangulation which is incident to two triangles. Therefore we obtain an element in $HG_{Eis}([4]_3)$.

It is natural to ask whether all elements in $HG_{Eis}([4]_3)$ can be obtained in this manner. The answer is affirmative. Start with a triangulation of desired type and unfold it to the plane accordingly by Alexandrov Unfolding Process. The polygon you get has the properties described before. Glue it as before to get the triangulation back.

It is also natural to ask whether two different polygons satisfying above properties give rise to different triangulations. In this case answer is not affirmative. To see this,
first observe that any such polygon is uniquely determined by it's vertex $\alpha$. Let $\delta=e^{\frac{2\pi \sqrt{-1}}{6}}$. If we change $\alpha$
with $\delta \alpha$, original polygon will be rotated in counter-clockwise direction by an angle of $\frac{2\pi}{6}$ around the origin. Therefore triangulation will not be changed. 
Hence we have a map
\begin{align}
\mathbf{Eis}/\langle \delta \rangle \rightarrow HG_{Eis}([4]_3)
\end{align}
and indeed, this map is also injective.

Observe that area of the polygon is proportional to the square-norm $\alpha \overline{\alpha}$, hence, square-norm of a lattice point gives number of triangles in the triangulation.
This case, $HG_{Eis}([4]_3)$, is also explained in \cite{shapes}.

We summarize the results of this section in the following theorem.
\begin{theorem}
There is a bijection
\begin{align}
\mathbf{Eis}/\langle \delta \rangle \equiv HG_{Eis}([4]_3)
\end{align}
such that square-norm of a lattice point gives number of triangles in corresponding triangulation.
\end{theorem}

\subsection{Shapes of Tetrahedra}
Let $C(\pi,\pi,\pi,\pi)$ \index{$C(\pi,\pi,\pi,\pi)$} be the set of cone metrics on sphere with four singular points of cone angle $\pi$, up to homotety and orientation preserving isometry.
The aim of this section is to describe this set.

Consider the following complex vector space
\begin{align}
H=\{(z_1,z_2):z_1,z_2 \in \mathbb{C}\},
\end{align}
with the Hermitian form
\begin{align}
\langle (z_1,z_2),(w_1,w_2)\rangle=\frac{\sqrt{-1}}{4}\{z_1\bar{w_2}-z_2\bar{w_1}\}.
\end{align}
If we regard an element $(z_1,z_2)$ as triangle in complex plane with vertices $0, z_1, z_2$, the square-norm of $(z_1,z_2)$
\begin{align*}
\frac{\sqrt{-1}}{4}\{z_1\bar{z_2}-z_2\bar{z_1}\}.
\end{align*}
gives signed area of the triangle, see Figure~\ref{ucgen}.
\definecolor{zzttqq}{rgb}{0.6,0.2,0}
\definecolor{qqqqff}{rgb}{0,0,1}
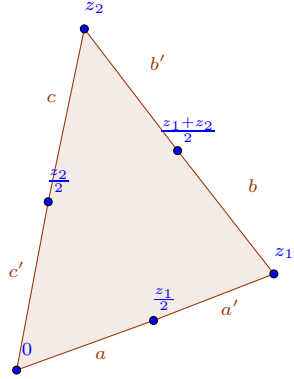
\begin{figure} \centering
\begin{tikzpicture}[line cap=round,line join=round,>=triangle 45,x=1.0cm,y=1.0cm]
\fill[color=zzttqq,fill=zzttqq,fill opacity=0.1] (1.66,0.2) -- (3.48,0.86) -- (5.08,1.48) -- (3.8,3.12) -- (2.56,4.74) -- (2.08,2.44) -- cycle;
\draw [color=zzttqq] (1.66,0.2)-- (3.48,0.86);
\draw [color=zzttqq] (3.48,0.86)-- (5.08,1.48);
\draw [color=zzttqq] (5.08,1.48)-- (3.8,3.12);
\draw [color=zzttqq] (3.8,3.12)-- (2.56,4.74);
\draw [color=zzttqq] (2.56,4.74)-- (2.08,2.44);
\draw [color=zzttqq] (2.08,2.44)-- (1.66,0.2);
\begin{scriptsize}
\draw [fill=qqqqff] (1.66,0.2) circle (1.5pt);
\draw[color=qqqqff] (1.8,0.48) node {$0$};
\draw [fill=qqqqff] (3.48,0.86) circle (1.5pt);
\draw[color=qqqqff] (3.62,1.14) node {$\frac{z_1}{2}$};
\draw [fill=qqqqff] (5.08,1.48) circle (1.5pt);
\draw[color=qqqqff] (5.22,1.76) node {$z_1$};
\draw [fill=qqqqff] (3.8,3.12) circle (1.5pt);
\draw[color=qqqqff] (3.94,3.4) node {$\frac{z_1+z_2}{2}$};
\draw [fill=qqqqff] (2.56,4.74) circle (1.5pt);
\draw[color=qqqqff] (2.7,5.02) node {$z_2$};
\draw [fill=qqqqff] (2.08,2.44) circle (1.5pt);
\draw[color=qqqqff] (2.22,2.72) node {$\frac{z_2}{2}$};
\draw[color=zzttqq] (2.78,0.4) node {$a$};
\draw[color=zzttqq] (4.5,1.04) node {$a'$};
\draw[color=zzttqq] (4.8,2.66) node {$b$};
\draw[color=zzttqq] (3.54,4.3) node {$b'$};
\draw[color=zzttqq] (2.12,3.82) node {$c$};
\draw[color=zzttqq] (1.66,1.54) node {$c'$};
\end{scriptsize}
\end{tikzpicture}
\caption{Cone Metric from an element in $H^+$}
\label{ucgen}
\end{figure}
Since there are both triangles of positive area and triangles of negative area, signature of this area Hermitian form is $(1,1)$. Let
\begin{align}
H^+=\{z\in H:\langle z,z \rangle >0\}.
\end{align}
be the positive part of $H$ with respect to area Hermitian form. $H^+$ consists of positively oriented triangles. There is a nice way to obtain cone metrics from these triangles.
Consider the triangle in Figure~\ref{ucgen} again. Glue the line segment $a$ with $a'$ and $b$ with $b'$. By this way we obtain a cone metric on sphere. It is clear that
angles at the vertices corresponding to $\frac{z_1}{2}, \frac{z_2}{2}, \frac{z_1+z_2}{2}$ are $\pi$. Also observe that the vertices $0, z_1, z_2$ come together to form a vertex
having angle $\pi$. Therefore we get an element in $C(\pi,\pi,\pi,\pi)$.
\definecolor{zzttqq}{rgb}{0.6,0.2,0}
\definecolor{qqqqff}{rgb}{0,0,1}
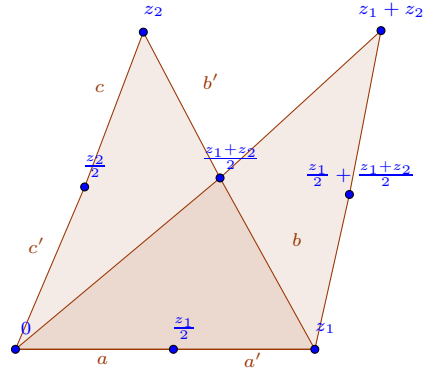
\begin{figure} \centering
\begin{tikzpicture}[line cap=round,line join=round,>=triangle 45,x=1.0cm,y=1.0cm]
\fill[color=zzttqq,fill=zzttqq,fill opacity=0.1] (1.12,-0.48) -- (3.22,-0.48) -- (5.1,-0.48) -- (3.84,1.8) -- (2.82,3.74) -- (2.04,1.68) -- cycle;
\fill[color=zzttqq,fill=zzttqq,fill opacity=0.1] (1.12,-0.48) -- (3.84,1.8) -- (5.98,3.76) -- (5.56,1.58) -- (5.1,-0.48) -- (3.22,-0.48) -- cycle;
\draw [color=zzttqq] (1.12,-0.48)-- (3.22,-0.48);
\draw [color=zzttqq] (3.22,-0.48)-- (5.1,-0.48);
\draw [color=zzttqq] (5.1,-0.48)-- (3.84,1.8);
\draw [color=zzttqq] (3.84,1.8)-- (2.82,3.74);
\draw [color=zzttqq] (2.82,3.74)-- (2.04,1.68);
\draw [color=zzttqq] (2.04,1.68)-- (1.12,-0.48);
\draw [color=zzttqq] (1.12,-0.48)-- (3.84,1.8);
\draw [color=zzttqq] (3.84,1.8)-- (5.98,3.76);
\draw [color=zzttqq] (5.98,3.76)-- (5.56,1.58);
\draw [color=zzttqq] (5.56,1.58)-- (5.1,-0.48);
\draw [color=zzttqq] (5.1,-0.48)-- (3.22,-0.48);
\draw [color=zzttqq] (3.22,-0.48)-- (1.12,-0.48);
\begin{scriptsize}
\draw [fill=qqqqff] (1.12,-0.48) circle (1.5pt);
\draw[color=qqqqff] (1.26,-0.2) node {$0$};
\draw [fill=qqqqff] (3.22,-0.48) circle (1.5pt);
\draw[color=qqqqff] (3.36,-0.2) node {$\frac{z_1}{2}$};
\draw [fill=qqqqff] (5.1,-0.48) circle (1.5pt);
\draw[color=qqqqff] (5.24,-0.2) node {$z_1$};
\draw [fill=qqqqff] (3.84,1.8) circle (1.5pt);
\draw[color=qqqqff] (3.98,2.08) node {$\frac{z_1+z_2}{2}$};
\draw [fill=qqqqff] (2.82,3.74) circle (1.5pt);
\draw[color=qqqqff] (2.96,4.02) node {$z_2$};
\draw [fill=qqqqff] (2.04,1.68) circle (1.5pt);
\draw[color=qqqqff] (2.18,1.96) node {$\frac{z_2}{2}$};
\draw[color=zzttqq] (2.28,-0.62) node {$a$};
\draw[color=zzttqq] (4.28,-0.62) node {$a'$};
\draw[color=zzttqq] (4.86,0.98) node {$b$};
\draw[color=zzttqq] (3.72,3.08) node {$b'$};
\draw[color=zzttqq] (2.24,3) node {$c$};
\draw[color=zzttqq] (1.4,0.9) node {$c'$};
\draw [fill=qqqqff] (5.98,3.76) circle (1.5pt);
\draw[color=qqqqff] (6.12,4.04) node {$z_1+z_2$};
\draw [fill=qqqqff] (5.56,1.58) circle (1.5pt);
\draw[color=qqqqff] (5.7,1.86) node {$\frac{z_1}{2}+\frac{z_1+z_2}{2}$};
\draw[color=zzttqq] (2.8,0.58) node {$$};
\draw[color=zzttqq] (5.24,2.72) node {$$};
\draw[color=zzttqq] (5.56,2.9) node {$$};
\draw[color=zzttqq] (5.12,0.78) node {$$};
\end{scriptsize}
\end{tikzpicture}
\caption{First Cutting Operation}
\label{cut1}
\end{figure}
Can every element in $C(\pi,\pi,\pi,\pi)$ be obtained from an element in $H^+$ by using the process above? Indeed, by Alexandrov Unfolding Process, we can cut-open an element in $C(\pi,\pi,\pi,\pi)$  to a polygon, actually a triangle, in $H^+$. We can glue edges of this triangle to get the cone metric back. Therefore we have a surjective map
\begin{align}
H^+\rightarrow C(\pi, \pi, \pi, \pi).
\end{align}
This map is far away from being injective. Let $\alpha \in \mathbb{C}$ be a complex number and $(z_1, z_2) \in H^+$. The triangle
\begin{align*}
\alpha(z_1,z_2)=(\alpha z_1, \alpha z_2)
\end{align*}
is obtained by rotating (around origin) and rescaling the triangle $ (z_1, z_2)$. Therefore triangles $(z_1, z_2)$ and $\alpha (z_1, z_2)$ give rise to the same element in
$C(\pi, \pi, \pi, \pi)$. Hence we have a map 
\begin{align}
\mathbb{P}H^+=\mathbb{H}=\mathbb{H}_{\mathbb{C}}^{1}\rightarrow C(\pi, \pi, \pi, \pi).
\end{align}
where $\mathbb{P}H^+$ is complex projectification of $H^+$ which is same as one dimensional complex hyperbolic space and 2 dimensional real hyperbolic space.

This map is not injective neither. Consider  Figure~\ref{cut1}. Given the triangle $(z_1, z_2)$ in $H^+$, we cut it through the line segment
$[0, \frac{z_1+z_2}{2}]$ and glue edges $b$ with $b'$ by a rotation of angle $\pi$ around $\frac{z_1+z_2}{2}$ to get the triangle $(z_1, z_1+z_2)$. Observe that the following elements gives the same cone metric:
$$
\begin{pmatrix}
  z_1 \\
  z_2
 \end{pmatrix}
\,\mbox{ and }\,
\begin{pmatrix}
  1 & 0 \\
  1 & 1
 \end{pmatrix}
 \begin{pmatrix}
  z_1 \\
  z_2
 \end{pmatrix}
 =
\begin{pmatrix}
  z_1 \\
  z_1+z_2
 \end{pmatrix}
 $$
\definecolor{zzttqq}{rgb}{0.6,0.2,0}
\definecolor{qqqqff}{rgb}{0,0,1}
\begin{figure} \centering
\begin{tikzpicture}[line cap=round,line join=round,>=triangle 45,x=1.0cm,y=1.0cm]
\fill[color=zzttqq,fill=zzttqq,fill opacity=0.1] (-0.12,0.84) -- (1.6,0.88) -- (3.32,0.92) -- (2.38,2.9) -- (1.68,4.44) -- (0.86,2.68) -- cycle;
\fill[color=zzttqq,fill=zzttqq,fill opacity=0.1] (1.68,4.44) -- (1.6,0.88) -- (1.52,-1.24) -- (0.76,-0.36) -- (-0.12,0.84) -- (0.86,2.68) -- cycle;
\fill[color=zzttqq,fill=zzttqq,fill opacity=0.1] (-0.12,0.84) -- (0.86,2.68) -- (1.68,4.44) -- (1.6,0.88) -- cycle;
\draw [color=zzttqq] (-0.12,0.84)-- (1.6,0.88);
\draw [color=zzttqq] (1.6,0.88)-- (3.32,0.92);
\draw [color=zzttqq] (3.32,0.92)-- (2.38,2.9);
\draw [color=zzttqq] (2.38,2.9)-- (1.68,4.44);
\draw [color=zzttqq] (1.68,4.44)-- (0.86,2.68);
\draw [color=zzttqq] (0.86,2.68)-- (-0.12,0.84);
\draw [color=zzttqq] (1.68,4.44)-- (1.6,0.88);
\draw [color=zzttqq] (1.6,0.88)-- (1.52,-1.24);
\draw [color=zzttqq] (1.52,-1.24)-- (0.76,-0.36);
\draw [color=zzttqq] (0.76,-0.36)-- (-0.12,0.84);
\draw [color=zzttqq] (-0.12,0.84)-- (0.86,2.68);
\draw [color=zzttqq] (0.86,2.68)-- (1.68,4.44);
\draw [color=zzttqq] (-0.12,0.84)-- (0.86,2.68);
\draw [color=zzttqq] (0.86,2.68)-- (1.68,4.44);
\draw [color=zzttqq] (1.68,4.44)-- (1.6,0.88);
\draw [color=zzttqq] (1.6,0.88)-- (-0.12,0.84);
\begin{scriptsize}
\draw [fill=qqqqff] (-0.12,0.84) circle (1.5pt);
\draw[color=qqqqff] (0.02,1.12) node {$0$};
\draw [fill=qqqqff] (1.6,0.88) circle (1.5pt);
\draw[color=qqqqff] (1.74,1.16) node {$\frac{z_1}{2}$};
\draw [fill=qqqqff] (3.32,0.92) circle (1.5pt);
\draw[color=qqqqff] (3.46,1.2) node {$z_1$};
\draw [fill=qqqqff] (2.38,2.9) circle (1.5pt);
\draw[color=qqqqff] (2.52,3.18) node {$\frac{z_1+z_2}{2}$};
\draw [fill=qqqqff] (1.68,4.44) circle (1.5pt);
\draw[color=qqqqff] (1.82,4.72) node {$z_2$};
\draw [fill=qqqqff] (0.86,2.68) circle (1.5pt);
\draw[color=qqqqff] (1,2.96) node {$\frac{z_2}{2}$};
\draw[color=zzttqq] (0.86,0.72) node {$a$};
\draw[color=zzttqq] (2.58,0.76) node {$a'$};
\draw[color=zzttqq] (3.24,2.22) node {$b$};
\draw[color=zzttqq] (2.44,3.98) node {$b'$};
\draw[color=zzttqq] (1.08,3.86) node {$c$};
\draw[color=zzttqq] (0.2,2.08) node {$c'$};
\draw [fill=qqqqff] (1.52,-1.24) circle (1.5pt);
\draw[color=qqqqff] (1.66,-0.96) node {$z_1-z_2$};
\draw [fill=qqqqff] (0.76,-0.36) circle (1.5pt);
\draw[color=qqqqff] (0.9,-0.08) node {$\frac{z_1-z_2}{2}$};
\draw[color=zzttqq] (1.36,0) node {$$};
\draw[color=zzttqq] (1.5,-0.42) node {$$};
\draw[color=zzttqq] (0.68,0.6) node {$$};
\end{scriptsize}
\end{tikzpicture} 
\caption{Second Cutting Operation}
\label{cut2}
\end{figure}
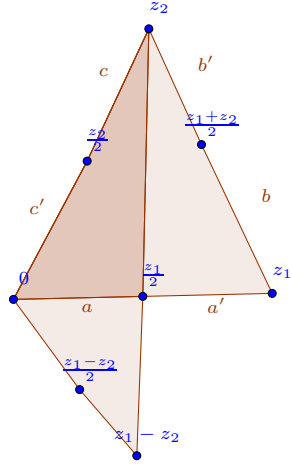
 Now consider Figure~\ref{cut2}. Cut the triangle $(z_1, z_2)$ from the line segment $[z_2, \frac{z_1}{2}]$ and glue $a$ with $a'$ by a rotation of angle $\pi$ around
$\frac{z_1}{2}$. You will get the triangle $(z_1-z_2, z_2)$ as in Figure~\ref{cut2}.

Observe that the following elements  gives the same cone metric. 
$$
\begin{pmatrix}
  z_1 \\
  z_2
 \end{pmatrix}
\,\mbox{  and }\, 
\begin{pmatrix}
  1 & -1 \\
  0 & 1
 \end{pmatrix}
\begin{pmatrix}
  z_1 \\
  z_2
 \end{pmatrix}
 =
\begin{pmatrix}
  z_1-z_2 \\
  z_2
 \end{pmatrix}
$$
The group generated by the matrices
 \begin{align}
  \begin{pmatrix}
  1 & 0 \\
  1 & 1
 \end{pmatrix},
 \begin{pmatrix}
  1 & -1 \\
  0 & 1
 \end{pmatrix}
\end{align}
 in $\mathrm{PSL}(2,\mathbb{R})$  is the modular group $\mathrm{PSL}(2,\mathbb{Z})$. Thus there is a well-defined map
 \begin{align}
  \mathbb{P}H^+/\mathrm{PSL}(2,\mathbb{Z})\rightarrow C(\pi, \pi, \pi, \pi)
   \end{align}
 We summarize results obtained in this section as follows:
    \begin{theorem}
    There is a map
\begin{align}
  \mathbb{P}H^+/\mathrm{PSL}(2,\mathbb{Z})\rightarrow C(\pi, \pi, \pi, \pi)
   \end{align}
   which is both injective and surjective.
\end{theorem}
 This bijection is not just a set theoretic bijection: one can naturally give, in some sense, complex structures to both    $\mathbb{P}H^+/\mathrm{PSL}(2,\mathbb{Z})$ and $C(\pi, \pi, \pi, \pi)$. The bijection above respects these structures. Also observe that above theorem means that $C(\pi, \pi, \pi, \pi)$ is nothing else than the modular orbifold.
\subsection{Back to triangulations}    
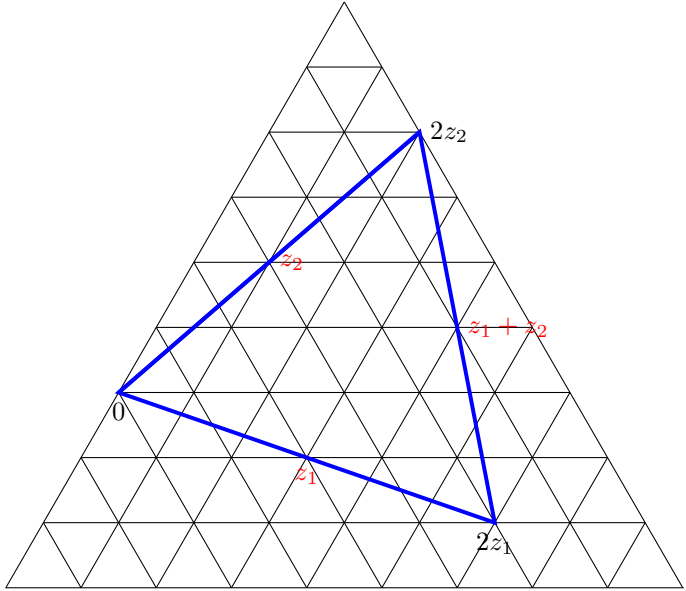
\begin{figure}
\newcommand*\rows{9}
 \centering

\begin{tikzpicture} 
 \coordinate (A) at  (1.5,{1.5*sqrt(3)})  ;
    \coordinate (C) at (6.5,{0.5*sqrt(3)});
     \coordinate (D) at (5.5,{3.5*sqrt(3)}) ;
               \foreach \row in {0, 1, ...,\rows} {
            \draw ($\row*(0.5, {0.5*sqrt(3)})$) -- ($(\rows,0)+\row*(-0.5, {0.5*sqrt(3)})$) ;    
            \draw ($\row*(1, 0)$) -- ($(\rows/2,{\rows/2*sqrt(3)})+\row*(0.5,{-0.5*sqrt(3)})$);
        \draw ($\row*(1, 0)$) -- ($(0,0)+\row*(0.5,{0.5*sqrt(3)})$);
    }
        \draw [blue,line width=1.5pt] (A) -- (C) node[,red,midway,below] {${z_1}$}    -- (D)  node[red,midway,right] {$z_1+z_2$}    -- cycle;  
  \draw [blue,line width=1.5pt] (A) --( A)  node[black,midway,below] {$0$} -- (D)  node[red,midway,right] {$z_2$} -- (D)  node[black,midway,right] {2$z_2$}  ;     
   
\draw [blue,line width=1.5pt] (C) --( C)  node[black,midway,below] {2$z_1$} ;

 \end{tikzpicture}
    \caption{Triangulations from Eisenstein Lattice} \label{E+}
    \end{figure}

 Set $\mathbb{E}^+=\mathbb{E}\cap H$, where 
\begin{align*}
\mathbb{E}=2\mathbf{Eis}\bigoplus 2\mathbf{Eis}=\{z=(2z_1,2z_2): z_i \in \mathbf{Eis} \} \subset H.
\end{align*}
Our next objective is to derive triangulations of type $[3]_4$ from $\mathbb{E}^+$.  Observe that the elements in $\mathbb{E}^+$  can be thought as positively oriented triangles whose vertices and midpoints of the edges are in $\mathbf{Eis}$. See Figure~\ref{E+}.
     There is a well-defined map
\begin{align}
\mathbb{E}^+\rightarrow HG_{Eis}([3]_4),
\end{align}
 given as follows. Take an element $(2 z_1, 2z_2)$ as in Figure~\ref{E+}. Glue the segment $[0,z_1]$ with $[z_1,2z_1]$ as we did before. Do the same for the segments
 $[2z_1, z_1+z_2]$, $[z_1+z_2, 2z_2]$ and $[2z_2,z_2]$, $[0,z_2]$. By this way, we get a triangulation of the sphere. Observe that the vertices incident to $z_1, z_2, z_1+z_2$
 are incident to 3 triangles. The vertices $0,2 z_1, 2z_2$  come together to form just one vertex of the triangulation which is also incident to 3 triangles. Therefore we obtain an element in $HG_{Eis}([3]_4)$.
 
 This map is surjective by Alexandrov Unfolding Process, one can cut-open a triangulation to obtain an element in $\mathbb{E}^+$ and glue this element appropriately 
 to get the initial triangulation back.
 
 This map is not injective. First of all, if $(2z_1, 2z_2) \in \mathbb{E}^+$ is a triangle, then multiplication by $\delta=e^{\frac{2\pi \sqrt{-1}}{6}}$ transforms this triangle 
 to $(2\delta z_1, 2\delta z_2) \in \mathbb{E}^+$ which is a triangle obtained by rotating the former triangle by an angle of $\frac{2\pi}{6}$ around origin. Therefore it does not change triangulation. Also cutting and gluing operations defined in the pervious section respect triangulations. Hence we have a map
 \begin{align}
\mathbb{E}^+/\langle \delta \rangle \times \mathrm{SL}(2,\mathbb{Z})\rightarrow HG_{Eis}([3]_4).
\end{align}
 This map is both injective and surjective.
 The area Hermitian form
 \begin{align}
\langle (z_1,z_2),(w_1,w_2)\rangle=\frac{\imath}{4}\{z_1\bar{w_2}-z_2\bar{w_1}\}.
\end{align} 
 defined in the previous section gives the area of the triangle considered. Therefore if we restrict our attention to $\mathbb{E}^+$, it gives us number of triangles in corresponding triangulation. Next theorem summarizes the results obtained in this section.  See also \cite{ismail}, \cite{ismailthesis}, \cite{shapes}.
 \begin{theorem}
There is a bijection
\begin{align}
\mathbb{E}^+/\langle \delta \rangle \times \mathrm{SL}(2,\mathbb{Z})\equiv  HG_{Eis}([3]_4),
\end{align}
such that the square-norm of each element gives number of triangles in the triangulation.
\end{theorem} 
 \subsection{Shapes of quadrangulations}\index{Shapes of quadrangulations}
Set $\mathbb{G}^+=\mathbb{G}\cap H$, where $\mathbb{G}$ is given by 
 \begin{align*}
 \mathbb{G}= 2\mathbb{Z}[\sqrt{-1}]\bigoplus 2\mathbb{Z}[\sqrt{-1}] 
 =\{ z=(2z_1, 2z_2): z_1, z_2 \in \mathbb{Z}[\sqrt{-1}]\}  \subset H
 \end{align*}
 We will obtain quadrangulations of type $[2]_4$ from $\mathbb{G}^+$. We consider quadrangulations as cone metrics by assuming that each quadrangle is unit  square. Observe that element in $\mathbb{G}^+$ can be regarded as triangles in complex plane whose vertices and midpoints of the edges are in the ring of Gaussian integers; $\mathbb{Z}[\sqrt{-1}]$. See Figure~\ref{G+}.
 
 \pagestyle{empty}
\definecolor{wwqqzz}{rgb}{0.4,0,0.6}
\definecolor{ffqqqq}{rgb}{1,0,0}
\definecolor{qqqqff}{rgb}{0,0,1}
\definecolor{cqcqcq}{rgb}{0.75,0.75,0.75}
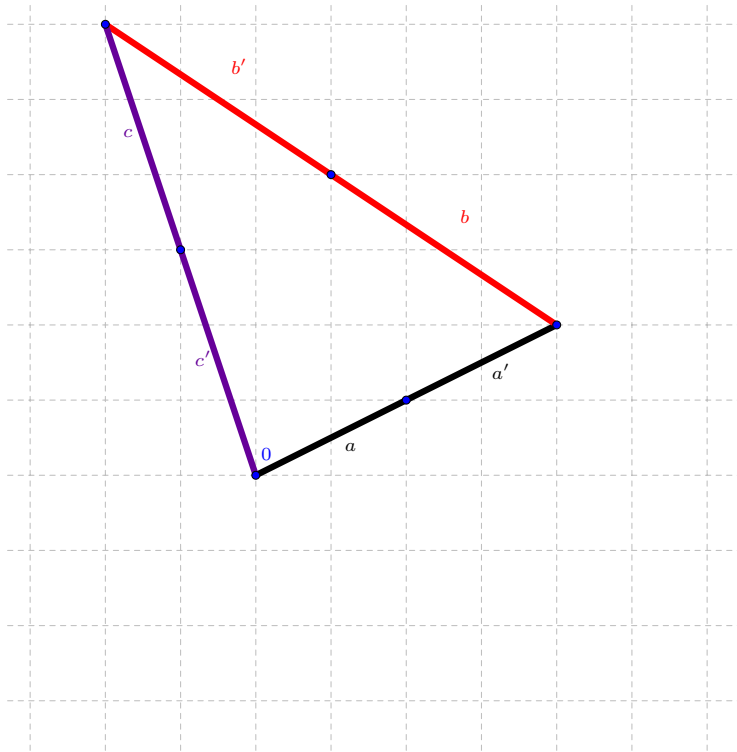
\begin{figure} \centering
\begin{tikzpicture}[line cap=round,line join=round,>=triangle 45,x=1.0cm,y=1.0cm]
\draw [color=cqcqcq,dash pattern=on 2pt off 2pt, xstep=1.0cm,ystep=1.0cm] (-4.3,-3.66) grid (5.48,6.3);
\draw [line width=2.4pt] (-1,0)-- (1,1);
\draw [line width=2.4pt] (1,1)-- (3,2);
\draw [line width=2.4pt,color=ffqqqq] (3,2)-- (0,4);
\draw [line width=2.4pt,color=ffqqqq] (0,4)-- (-3,6);
\draw [line width=2.4pt,color=wwqqzz] (-3,6)-- (-2,3);
\draw [line width=2.4pt,color=wwqqzz] (-2,3)-- (-1,0);
\begin{scriptsize}
\draw [fill=qqqqff] (-1,0) circle (1.5pt);
\draw[color=qqqqff] (-0.86,0.28) node {$0$};
\draw [fill=qqqqff] (1,1) circle (1.5pt);
\draw[color=black] (0.26,0.38) node {$a$};
\draw [fill=qqqqff] (3,2) circle (1.5pt);
\draw[color=black] (2.26,1.38) node {$a'$};
\draw [fill=qqqqff] (0,4) circle (1.5pt);
\draw[color=ffqqqq] (1.78,3.44) node {$b$};
\draw [fill=qqqqff] (-3,6) circle (1.5pt);
\draw[color=ffqqqq] (-1.22,5.44) node {$b'$};
\draw [fill=qqqqff] (-2,3) circle (1.5pt);
\draw[color=wwqqzz] (-2.7,4.56) node {$c$};
\draw[color=wwqqzz] (-1.7,1.56) node {$c'$};
\end{scriptsize}
\end{tikzpicture}
\caption{Quadrangulations from Gaussian Lattice}
\label{G+}
\end{figure}
Gluing process explained before provides us a map 
\begin{align}
\mathbb{G}^+ \rightarrow HG_{Gauss}([2]_4).
\end{align}
This map is surjective by Alexandrov Unfolding Process, but it is not injective. Multiplication of a lattice element by $\sqrt{-1}$ just rotates the triangle by an angle of 
$\frac{\pi}{2}$ around origin; thus it respects quadrangulation. Also cutting and gluing operations defined before respect quadrangulation. We have a map
 \begin{align}
\mathbb{G}^+/\langle \sqrt{-1} \rangle \times \mathrm{SL}(2,\mathbb{Z})\rightarrow HG_{Gauss}([2])_4 .
\end{align}
 which is both injective and surjective.
 
 Observe that area Hermitian form give the number of quadrangles in corresponding quadrangulation. Following theorem summarizes the results obtained in this section.
 \begin{theorem}
There is a bijection
\begin{align}
\mathbb{G}^+/\langle \sqrt{-1} \rangle \times \mathrm{SL}(2,\mathbb{Z})\equiv  HG_{Gauss}([2]_4)
\end{align}
such that square-norm of each element gives number of quadrangles in the quadrangulation.
\end{theorem} 

\section{Beyond  Hypergeometric}
It must be possible to extend the classification of results of triangulations of non-negative curvature to more general triangulations (same for the quadrangulations).  To achieve this, we need the right conditions to control the curvature. Suggestions: ``just one point of negative curvature above infinity", or ``just one point of negative curvature above infinity, whose curvature is bounded below by $\kappa$ ", or ``just one point of fixed curvature $\kappa$ above infinity" (in each case, the points of non-negative curvature are arbitrary). 
We may also allow for a fixed number of points with controlled negative curvature. These relaxed conditions may bring in non-discrete groups into the picture, the signatures of the Hermitian forms will change, complex hyperbolic structure will decay, and there is a possibility that the parameter spaces will brake up into disconnected components. On the other hand, the relaxed conditions may lead to the discovery of other arithmetic and non-arithmetic discrete groups acting on some symmetric or non-symmetric spaces, e.g. ``complex deSitter spaces".

If we further relax the control of the points of negative curvature by simply requiring that it be bounded globally from below, then things will totally go out of control. 
Indeed it is easy to illustrate how wild things may become in terms of quadrangulations. 
Consider a big cube as in Figure \ref {bigcube}, with $6n^2$ quadrangles. 
Its surface is a hypergeometric sphere quadrangulation (there are 8 points of positive curvature). Now imagine that this cube is made of $n^3$ smaller cubes. Imagine that you are a sculpture. Then by removing smaller cubes you may carve out any three-dimensional figure with galleries inside, and the curvature will remain bounded below by $-\pi/2-$ -which is already the greatest negative value that the curvature may attain in this case. In fact you may decide to glue little cubes to form self-overlaps of the 3-d figure. 

It is much harder to describe the situation as the curvature goes deeper, since the shapes become non-embeddable locally in this case.
So, it seems that abolishing all restrictions on the curvature (including the condition of being bounded from below) do not lead to a well-posed problem, neither.

\bigskip
\noindent\ITEM It might be appropriate to conclude this text with an apology:
The term ``hypergeometric curve" is used in the literature to refer to some {\it families} of cyclic branched coverings of
the projective line. Here this term refers to certain {\it rigid} (arithmetic) curves, which can be described by some special dessins
(equivalently by triangulations, origamis, quilts, etc). Since this terminology seems to unify the rich vocabulary surrounding the hypergeometric
phenomena, we could not resist the temptation to call these curves {\it hypergeometric}.

\medskip
\noindent {\bf Acknowledgements.} 
We are thankful to Athanase Papadopoulos for inviting us to publish in this volume. 
Both authors were funded by the grant TUBITAK-110T690.
The first named is author was funded by a Galatasaray University Research Grant and 
the grant TUBITAK-114R073.

\label{references}

\input{HGHTbiblio}

\newpage
\section* {APPENDIX 1}

 \begin{center}
 \scalebox{0.85}{
\begin{tabular}{|c|ccccc|c|c|c|c|c|}
\hline
dim & $k_1$ & $k_2$ & $k_3$ & $k_4$ & $k_5$ &\!deg\!&\!\!Compct?\!\!&\!Number\!&\!Pure?\!&\!ar?\!\\\hline
9&0&0&0&0&12&2&N&10&I&AR\\\hline
8&0&0&0&1&10&2&N&11&I&AR\\\hline
7&0&0&1&0&9&2&N&12&I&AR\\\hline
7&0&0&0&2&8&2&N&13&I&AR\\\hline
6&0&1&0&0&8&2&N&14&I&AR\\\hline
6&0&0&1&1&7&2&N&15&I&AR\\\hline
5&1&0&0&0&7&2&N&16&I&AR\\\hline
6&0&0&0&3&6&2&N&17&I&AR\\\hline
5&0&1&0&1&6&2&N&18&I&AR\\\hline
5&0&0&2&0&6&2&N&19&I&AR\\\hline
5&0&0&1&2&5&2&N&20&I&AR\\\hline
4&1&0&0&1&5&2&N&22&I&AR\\\hline
4&0&1&1&0&5&2&N&23&I&AR\\\hline
5&0&0&0&4&4&2&N&24&I&AR\\\hline
4&0&0&2&1&4&2&N&25&I&AR\\\hline
3&1&0&1&0&4&2&N&26&I&AR\\\hline
3&0&2&0&0&4&2&N&27&I&AR\\\hline
4&0&0&1&3&3&2&N&28&I&AR\\\hline
3&1&0&0&2&3&2&N&29&I&AR\\\hline
3&0&1&1&1&3&2&N&30&I&AR\\\hline
3&0&0&3&0&3&2&N&31&I&AR\\\hline
3&0&0&0&6&0&2&N&1&P&AR\\\hline
2&0&1&0&4&0&2&N&2&P&AR\\\hline
2&1&1&0&0&3&2&N&32&I&AR\\\hline
4&0&0&0&5&2&2&N&33&I&AR\\\hline
4&0&2&0&3&2&2&N&34&I&AR\\\hline
3&0&0&2&2&2&2&N&35&I&AR\\\hline
2&1&0&1&1&2&2&N&36&I&AR\\\hline
2&0&2&0&1&2&2&N&37&I&AR\\\hline
2&1&0&2&0&2&2&N&38&I&AR\\\hline
3&0&0&1&4&1&2&N&39&P&AR\\\hline
2&1&0&0&3&1&2&N&40&P&AR\\\hline
2&0&1&1&2&1&2&N&41&P&AR\\\hline
2&0&0&3&1&1&2&N&42&P&AR\\\hline
2&0&0&2&3&0&2&N&43&P&AR\\\hline
1&1&0&1&2&0&-&N&-&&AR\\\hline
1&1&0&2&0&1&-&N&-&&AR\\\hline
1&1&1&0&1&1&-&N&-&&AR\\\hline
1&0&1&2&1&0&-&N&-&&AR\\\hline
1&0&2&0&2&0&-&N&-&&AR\\\hline
1&0&2&1&0&1&-&N&-&&AR\\\hline
1&0&0&4&0&0&-&N&-&&AR\\\hline
0&1&1&1&0&0&-&self&-&&\\\hline
0&2&0&0&0&1&-&-&-&&\\\hline
0&2&0&0&1&0&-&-&-&&\\\hline
0&0&3&0&0&0&-&self&-&&\\\hline
\end{tabular}
}
\end{center}

\newpage
\section*{APPENDIX 2}
\begin{center}
\scalebox{0.85}{
\begin{tabular}{|c|cc|ccccc|}

\hline
dim & $m_1$ & $n_1$ & $k_1$ & $k_2$ & $k_3$ & $k_4$ & $k_5$\\\hline
*6&1&0&  0&0&0&0&9\\\hline
*5&1&0&  0&0&0&1&7\\\hline
*4&1&0&  0&0&1&0&6\\\hline
*4&1&0&  0&0&0&2&5\\\hline
*3&1&0&  0&1&0&0&5\\\hline
*3&1&0&  0&0&1&1&4\\\hline
*2&1&0&  1&0&0&0&4\\\hline
*3&1&0&  0&0&0&3&3\\\hline
*2&1&0&  0&0&2&0&3\\\hline
*2&1&0&  0&1&0&1&3\\\hline
*2&1&0&  0&0&1&2&2\\\hline
1&1&0&  0&1&1&0&2\\\hline
1&1&0&  1&0&0&1&2\\\hline
*2&1&0&  0&0&0&4&1\\\hline
1&1&0&  0&1&0&2&1\\\hline
0&1&0&  0&2&0&0&1\\\hline
0&1&0&  1&0&1&0&1\\\hline
1&1&0&  0&0&1&3&0\\\hline
0&1&0&  1&0&0&2&0\\\hline
0&1&0&  0&1&1&1&0\\\hline
-1&1&0&  1&1&0&0&0\\\hline
*3&2&0&  0&0&0&0&6\\\hline
*2&2&0&  0&0&0&1&4\\\hline
1&2&0&  0&0&1&0&3\\\hline
1&2&0&  0&0&0&2&2\\\hline
0&2&0&  0&1&0&0&2\\\hline
-1&2&0&  1&0&0&0&1\\\hline
0&2&0&  0&1&1&0&1\\\hline
0&2&0&  0&0&0&3&0\\\hline
-1&2&0&  0&0&2&0&0\\\hline
0&3&0&  0&0&0&0&3\\\hline
-1&3&0&  0&0&0&1&1\\\hline
-2&3&0&  0&0&1&0&0\\\hline
\end{tabular}
}
\scalebox{0.85}{
\begin{tabular}{|c|cc|ccccc|}
\hline
dim & $m_1$ & $n_1$ & $k_1$ & $k_2$ & $k_3$ & $k_4$ & $k_5$\\\hline
*5&0&1&  0&0&0&0&8\\\hline
*4&0&1&  0&0&0&1&6\\\hline
*3&0&1&  0&0&1&0&5\\\hline
*3&0&1&  0&0&0&2&4\\\hline
*2&0&1&  0&1&0&0&4\\\hline
*2&0&1&  0&0&1&1&3\\\hline
1&0&1&  1&0&0&0&3\\\hline
*2&0&1&  0&0&0&3&2\\\hline
1&0&1&  0&0&2&0&2\\\hline
1&0&1&  0&1&0&1&2\\\hline
1&0&1&  0&0&1&2&1\\\hline
0&0&1&  0&1&1&0&1\\\hline
0&0&1&  1&0&0&1&1\\\hline
1&0&1&  0&0&0&4&0\\\hline
0&0&1&  0&1&0&2&0\\\hline
-1&0&1&  0&2&0&0&0\\\hline
-1&0&1&  1&0&1&0&0\\\hline
1&0&2&  0&0&0&0&4\\\hline
1&0&2&  0&0&1&1&2\\\hline
0&0&2&  1&0&1&0&1\\\hline
-1&0&2&  0&0&0&2&0\\\hline
-2&0&2&  0&1&0&0&0\\\hline
-3&0&3&  0&0&0&0&0\\\hline
*2&1&1&  0&0&0&0&5\\\hline
1&1&1&  0&0&0&1&3\\\hline
0&1&1&  0&0&1&0&2\\\hline
-1&1&1&  0&1&0&0&1\\\hline
-2&1&1&  1&0&0&0&0\\\hline
-1&2&1&  0&0&0&0&2\\\hline
-2&2&1&  0&0&0&1&0\\\hline
-2&1&2&  0&0&0&0&1\\\hline
-3&4&0&  0&0&0&0&0\\\hline
&&&&&&&\\\hline
\end{tabular}
}
\end{center}

\newpage 
\section*{APPENDIX 3}

\begin{center}
\begin{tabular}{|cc|}
\hline
($m_1$, $n_1$, $n_2$)=(1,0,0) or (0,0,1)
&
\begin{tabular}{|c|ccc|}
dim &$k_1$ & $k_2$ & $k_3$\\\hline
3&0&0&6\\\hline
2&0&1&4\\\hline
2&1&0&3\\\hline
1&0&2&2\\\hline
0&1&1&1\\\hline
0&0&3&0\\\hline
-1&2&0&0
\end{tabular}\hspace{-4.2mm}\\\hline
($m_1$, $n_1$, $n_2$)=(0,1,0)
&
\begin{tabular}{|c|ccc|}
dim &$k_1$ & $k_2$ & $k_3$\\\hline
2&0&0&5\\\hline
1&0&1&3\\\hline
0&1&0&2\\\hline
0&0&2&1\\\hline
-1&1&1&0
\end{tabular}\hspace{-4.2mm} \\\hline
($m_1$, $n_1$, $n_2$)=(2,0,0), (0,0,2) or (1,0,1)
&
\begin{tabular}{|c|ccc|}
dim &$k_1$ & $k_2$ & $k_3$\\\hline
1&0&0&4\\\hline
0&0&1&2\\\hline
-1&1&0&1\\\hline
-1&0&2&0\\\hline
0&0&0&3\\\hline
-1&0&1&1\\\hline
-2&1&0&0
\end{tabular}\hspace{-4.2mm} \\\hline
($m_1$, $n_1$, $n_2$)=(1,1,0), (0,1,1)
&
\begin{tabular}{|c|ccc|}
dim &$k_1$ & $k_2$ & $k_3$\\\hline
0&0&0&3\\\hline
-1&0&1&1\\\hline
-2&1&0&0
\end{tabular}\hspace{-4.2mm} \\\hline
($m_1$, $n_1$, $n_2$)=(0,2,0)
&
\begin{tabular}{|c|ccc|}
dim &$k_1$ & $k_2$ & $k_3$\\\hline
-1&0&0&2\\\hline
-2&0&1&0
\end{tabular}\hspace{-4.2mm} \\\hline
($m_1$, $n_1$, $n_2$)=(2,1,0) or (0,1,2)
&
\begin{tabular}{|c|ccc|}
dim &$k_1$ & $k_2$ & $k_3$\\\hline
-2&0&0&1
\end{tabular}\hspace{-4.2mm} \\\hline
($m_1$, $n_1$, $n_2$)=(1,2,0), (0,2,1), (4,0,0), (0,0,4) 
&
\begin{tabular}{|c|ccc|}
dim &$k_1$ & $k_2$ & $k_3$\\\hline
-3&0&0&0
\end{tabular}\hspace{-4.2mm} \\\hline
\end{tabular}
\end{center}

\end{document}